\documentclass[a4paper,11pt]{article}
\usepackage{MyStyle}
\usepackage{etrees}

\allowdisplaybreaks

\newtheorem{theorem}{Theorem}[section]
\newtheorem{definition}[theorem]{Definition}
\newtheorem{proposition}[theorem]{Proposition}
\newtheorem{lemma}[theorem]{Lemma}
\newtheorem{remark}[theorem]{Remark}
\newtheorem*{remark*}{Remark}

\newtheorem{ass}[theorem]{Assumption}
\newtheorem*{notation*}{Notation}

\newtheorem*{ex*}{Example}

\newtheorem*{exs*}{Examples}

\newtheorem*{app*}{Application}

\def\ts{\thinspace}

\newcommand{\lpar}{\bigg(}
\newcommand{\rpar}{\bigg)}
\newcommand{\lcro}{\bigg[}
\newcommand{\rcro}{\bigg]}

\title{
Exotic aromatic B-series for the study of long time integrators for a class of ergodic SDEs
}

\author{ 
Adrien Laurent\textsuperscript{1,2} and Gilles Vilmart\textsuperscript{2}
}

\begin{document}
\footnotetext[1]{
Ecole Normale Sup\'erieure de Rennes and Universit\'e de Rennes 1,  
Campus de Ker Lann, avenue Robert Schumann, F-35170 Bruz, France, Adrien.Laurent@unige.ch}
\footnotetext[2]{
Universit\'e de Gen\`eve, Section de math\'ematiques, 2-4 rue du Li\`evre, CP 64, CH-1211 Gen\`eve 4, Switzerland, Gilles.Vilmart@unige.ch}

\maketitle

\begin{abstract}

We introduce a new algebraic framework based on a modification (called exotic) of aromatic Butcher-series for the systematic study of the accuracy of numerical integrators for the invariant measure of a class of ergodic stochastic differential equations (SDEs) with additive noise. The proposed analysis covers Runge-Kutta type schemes including the cases of partitioned methods and postprocessed methods. We also show that the introduced exotic aromatic B-series satisfy an isometric equivariance property.

\smallskip

 \end{abstract}

\section{Introduction}
We consider a class of stochastic systems of differential equations of the form
\begin{equation} \label{equation:sde}
dX(t) = f(X(t)) dt + \sigma dW(t),
\end{equation}
where $X(t)\in\R^d$ is the solution with initial condition $X_0$ assumed deterministic for simplicity, 
the vector field $f:\R^d\rightarrow \R^d$ is assumed smooth and globally Lipschitz, 
$\sigma>0$ is a constant,
and $W(t)$ is
a standard $d$-dimensional Wiener process fulfilling the usual assumptions.

We say that problem \eqref{equation:sde} is ergodic if it has a unique invariant measure $\mu$ satisfying for all deterministic initial conditions $X_0$ and all smooth test functions $\phi$,
\begin{equation} \label{equation:ergodicityT}
\lim_{T\rightarrow \infty} \frac1{T} \int_{0}^T\phi(X(s))ds = \int_{\R^d} \phi(x) d\mu(x), \qquad \mbox{almost surely}.
\end{equation}
Under appropriate smoothness and growth assumptions on the vector field $f$, the above ergodicity property is automatically satisfied, and in addition one has in general the following exponential convergence for all initial conditions $X_0$ and all appropriate test functions $\phi$,
\begin{equation} \label{equation:ergodicityexp}
\left|\E[\phi(X(t))] - \int_{\R^d} \phi(x) d\mu(x)\right| \leq Ce^{-\lambda t}
\end{equation}
for all $t>0$ where $C=C(\phi,X_0)$ and $\lambda>0$ are independent of time $t>0$.
We refer to \cite{Bally96tlo,Mattingly10con,Debussche12wbe,Kopec15wbea,Kopec15wbeb} for details. 
A special case of the class of problems \eqref{equation:sde} is the well-studied overdamped Langevin equation, also called Brownian dynamics,\footnote{Up to a time transformation, one could for simplicity fix the value of $\sigma$ (e.g. \ts $\sigma=1$ or $\sqrt{2}$), however we choose not to fix it to better distinguish the various order conditions.}
\begin{equation} \label{equation:Langevin}
dX(t) = -\nabla V(X(t)) dt + \sigma dW(t)
\end{equation}
where $V:\R^d\rightarrow \R$ is a smooth potential,  $f=-\nabla V$ and $\sigma>0$.
Under appropriate growth conditions on the potential $V$ (growing at least quadratically), \eqref{equation:Langevin} is ergodic and \eqref{equation:ergodicityT} and \eqref{equation:ergodicityexp} hold with the invariant measure $d\mu(x) = \rho_\infty(x) dx$ where its density $\rho_\infty$ with respect to the Lebesgue measure is given by the Gibbs distribution
\begin{equation} \label{eq:rho}
\rho_\infty(x) = Z e^{-\frac{2}{\sigma^{2}} V(x)}
\end{equation}
where $Z$ is a normalization constant such that $\int_{\R^d} \rho_\infty(x) dx =1$.
Computing integrals with respect to the invariant measure $\mu$ in high dimensions is in general very costly using a deterministic quadrature rule, and one can take advantage of the above ergodicity property to compute such integrals using numerical approximations of \eqref{equation:sde}, the exact solution of \eqref{equation:sde} not being available in general.
Let us mention that a natural way to sample without any bias from the invariant measure with the Gibbs density \eqref{eq:rho} is to apply Markov Chain Monte-Carlo methods, in particular
Metropolis–Hastings algorithms, see for instance the survey \cite{SanzSerna14mcm}. However, as highlighted in \cite{Milstein07cel}, ``the Metropolis–Hastings algorithm is rather expensive due to the need of accept/reject steps and does not admit the use of
powerful weak methods'', which can be combined with the methodology of ``rejecting exploding  trajectories'' \cite{Milstein07cel}. We thus focus in this paper on classes of weak methods and
consider a one step numerical integrator for the
approximation of \eqref{equation:sde} at time $t_n=nh$ of the form 
\begin{equation} \label{equation:defPsi}
X_{n+1} = \Psi(X_n,h,\xi_n),
\end{equation}
where $h$ is a fixed time step size and $\xi_n$ are independent random vectors.
We say that the numerical scheme has local weak order $r$ if the weak error after one step satisfies
\begin{equation}
\label{equation:weak_condition}
|\E[\phi(X_1)|X_0=x] - \E[\phi(X(h))|X(0)=x]| \leq C h^{r+1},
\end{equation}
where $C=C(\phi,x)$ is independent of $h$ assumed small enough and $\phi$ is a test function.
Note that under appropriate assumptions on the numerical scheme (to achieve in particular bounded numerical moments along time), one can in general deduce a global weak order $r$, $
|\E(\phi(X_n)) - \E(\phi(X(t_n)))| \leq C h^{r}
$, as shown in \cite{Milstein85wao} (see \cite[Chap.\ts 2.2]{Milstein04snf}).
The numerical method is called ergodic if it has a unique invariant probability law $\mu^h$ with finite moments of any order and
$$
\lim_{N\rightarrow \infty} \frac1{N+1} \sum_{n=0}^N \phi(X_n) =
\int_{\R^d} \phi(x) d\mu^h(x), \qquad \mbox{almost surely},
$$
for all deterministic initial conditions $X_0=x$ and all test functions $\phi$.
We say that the numerical method has order $p$ with respect to the invariant measure of \eqref{equation:sde} if
\begin{equation} \label{equation:order}
\left|\int_{\R^d} \phi(x) d\mu^h(x) - \int_{\R^d} \phi(x) d\mu(x)\right| \leq Ch^p,
\end{equation}
where $C$ is independent of $h$ assumed small enough.
Under appropriate assumptions on the numerical scheme (see for instance \cite{Debussche12wbe}), one also obtains
the following exponential estimate similar to \eqref{equation:ergodicityexp} (possibly with a different constant $\lambda>0$),
\begin{equation} \label{equation:orderexp}
\left|\E[\phi(X_n)] - \int_{\R^d} \phi(x) d\mu(x)\right| \leq Ke^{-\lambda t_n} + Ch^p,
\end{equation}
where $K=K(\phi,x)$, $C=C(\phi)$ are independent of $n$ and $h$ assumed small enough.
A simple way to achieve high order $p$ for the invariant measure is to consider a numerical scheme with high standard weak order $r$, and it is known for large classes of SDEs that $p\geq r$, see in particular \cite{Mattingly02efs} in the context of locally Lipschitz vector fields with multiplicative noise.
Note analogously that the strong order $q$ of convergence, which corresponds to the numerical approximation of individual trajectories of \eqref{equation:sde}, is in general lower than or equal to the weak order $r$ of convergence. 
There are interestingly many schemes in the literature for which $p>r$ and a high order $p$ for the invariant measure is obtained, while the standard weak order of accuracy remains low, typically of order $r=1$, i.e.\ts the scheme is consistent in the weak convergence sense. This is the case in particular for the Langevin equation \cite{BouRabee10lra,Leimkuhler13rco,Leimkuhler16tco,Abdulle15lta}. 
In \cite{Abdulle14hon,Abdulle15lta}, a methodology for the analysis and design of high order integrators for the invariant measure is introduced and serves as a crucial ingredient in this paper.
The approach combines the usual Talay-Tubaro methodology \cite{Talay90eot} and recent developments of the theory of backward error analysis and modified differential equations in the stochastic context \cite{Zygalakis11ote,Abdulle12hwo,Debussche12wbe,Kopec15wbea,Kopec15wbeb}, a major tool in the area of deterministic geometric numerical integration \cite{Hairer06gni}. 
In \cite{Vilmart15pif} for finite dimensions and in \cite{Brehier16hoi} in the context of parabolic stochastic partial differential equations, this approach is combined with the idea of processing from Butcher \cite{Butcher69teo}, to design efficient postprocessed integrators with high order for the invariant measure at a negligible overcost compared to standard low order schemes.
The postprocessor methodology is extended in \cite{Abdulle17oes} for a class of explicit stabilized schemes of order two for the invariant measure and with optimally large stability domains.

The aim of this paper is to provide a unified algebraic framework based on aromatic trees and B-series, with a set of trees independent of the dimension $d$ of the problem, for the systematic study of the order conditions for the invariant measure of a class of numerical integrators that includes Runge-Kutta type schemes for problems of the form \eqref{equation:sde}. 
We show that the new framework permits to recover some schemes and simplify the calculations in \cite{Abdulle14hon} and for postprocessed integrators in \cite {Vilmart15pif,Brehier16hoi,Abdulle17oes}.
Analogously to \cite{Rossler06rta} (we study here the additive noise case), we consider in this paper Runge-Kutta methods of the form\footnote{%
Note that the internal stages $Y_i$ depend on $n$, but this dependence is omitted for brevity of the notation.}
\begin{equation}
\label{equation:defRK}
\begin{array}{rl}
  Y_i&=X_n+h\sum \limits_{j=1}^s a_{ij}f(Y_j) +\sum \limits_{k=1}^l d_i^{(k)} \sigma \sqrt{h}\xi_n^{(k)} , \qquad i = 1, \dots ,s, \\
  X_{n+1}&=X_n+h\sum \limits_{i=1}^s b_i f(Y_i) +\sigma\sqrt{h}\xi_n^{(1)},
\end{array}
\end{equation}
where $a_{ij},b_i,d_i^{(k)}$ are the coefficients defining the Runge-Kutta scheme, and $\xi_n^{(k)}\sim\NN(0,I_d)$ are independent Gaussian random vectors.
We highlight once again that we focus in this paper on the high order $p$ of accuracy for the invariant measure, while the order $r$ of accuracy in the weak sense can remain low (typically $r=1$).
The analysis in this paper applies to the class of methods \eqref{equation:defRK} for any number $l$ of random vectors in the internal stages. However we shall often consider $l=1$ random vector per internal stage\footnote{In this case, we denote $d=d^{(1)}$ and $\xi_n=\xi_n^{(1)}$.}, which is sufficient to achieve order $p=2$ or $3$ for the invariant measure.
In particular, we shall consider the $\theta$-method as an illustrative example in Sections \ref{section:using_eat} and \ref{section:construction_high_order_integrators} and recover known results on its accuracy. It is defined for $\theta$ fixed as
\begin{equation}
\label{equation:theta_method}
X_{n+1}=X_n+ h(1-\theta)f(X_n) +h\theta f(X_{n+1})+\sigma \sqrt{h}\xi_n,
\end{equation}
For $\theta=0$, we get the explicit Euler-Maruyama method while the scheme is implicit for $\theta \neq 0$.  It can be put in Runge-Kutta form \eqref{equation:defRK} for $l=1$ with the following coefficients.
$$\begin{array}
{c|c|c}
c & A & d\\
\hline
  & b & 
\end{array}
=
\begin{array}
{c|cc|c}
0 & 0 & 0 & 0\\
1 & 1-\theta & \theta & 1\\
\hline
  & 1-\theta & \theta &  
\end{array}$$


The usage of trees and B-series\footnote{%
originally named Butcher-series} is known as a powerful standard tool for the numerical analysis of differential equations. 
B-series where introduced by Hairer and Wanner in \cite{Hairer74otb} based on the work of Butcher \cite{Butcher72aat}, and are now exposed in many articles and books \cite{Hairer06gni,Butcher16nmf}, see also the presentation in \cite{Chartier10aso,SanzSerna15fsa}.
In the last decades, several works extended trees and B-series to the stochastic context, we mention in particular Burrage and Burrage \cite{Burrage96hso} and Komori, Mitsui and Sugiura \cite{Komori97rta} who first introduced stochastic trees and B-series for studying the order conditions of strong convergence of SDEs, and \cite{Burrage00oco, Rossler04ste, Rossler06rkm, Rossler06rta, Debrabant08bsa, Rossler10saw, Debrabant11cos} for the design and analysis of high order weak and strong integrators on a finite time interval. Tree series were also used to describe schemes preserving quadratic invariants \cite{Anmarkrud17ocf}.
The recent work \cite{Alamo16atf} also studies algebraically strong and weak errors, but instead of using tree series, it uses word series because these are well suited in the context of splitting stochastic integrators.

In this paper, we focus on the long time accuracy of numerical integrators and derive in a systematic manner the order conditions for sampling the invariant measure of an ergodic system of the form \eqref{equation:sde}. Additionally, in Section \ref{Section:non_reversible_perturbation}, we allow the inclusion of non-reversible perturbation as in \cite{Lelievre13onr,Duncan16vru}.
The proposed algebraic framework relies on aromatic B-series, a generalisation of B-series introduced in \cite{MuntheKaas16abs, McLachlan16bsm} (see also the presentation in \cite{Bogfjellmo15aso}) to characterize all the schemes that are affine equivariant, i.e.\ts that behave transparently with respect to an affine change of coordinates.
These aromatic B-series rely on aromatic trees, which were first introduced in \cite{Chartier07pfi} to represent the divergence of B-series in the context of deterministic value or first integral preserving ordinary differential equations.

This paper is organized as follows.
In Section \ref{section:preliminaries}, we described the general setting and assumptions needed in our analysis.
In Section \ref{section:exotic_aromatic_forests}, we introduce a new generalization of B-series, called ``exotic aromatic B-series'' by considering an additional new type of edge called ``liana'' compared to standard aromatic B-series.
We also show that these new exotic aromatic B-series satisfy an isometric equivariance property.
In Section \ref{section:using_eat}, we explain how this new algebraic framework applies for the long time accuracy analysis of stochastic integrators for ergodic problems.
In Section \ref{section:construction_high_order_integrators}, we derive order conditions for integrators expandable as aromatic B-series methods, with special emphasis on Runge-Kutta type integrators and post-processed integrators.
In particular, we show that the orders 2 and 3 for the invariant measure of Brownian dynamics \eqref{equation:Langevin} yield respectively 2 and 6 order conditions (see Table \ref{table:RK_order_conditions}), compared to the 3 and 10 more restrictive conditions for the standard weak order of convergence (see Table \ref{table:RK_weak_order_conditions}).

\section{Preliminaries}
\label{section:preliminaries}
We first state the following smoothness and growth Assumptions \ref{Assumption:growth} and \ref{Assumption:f_gradient} on the vector field of \eqref{equation:sde} which automatically yield that \eqref{equation:sde} satifies the
ergodicity properties \eqref{equation:ergodicityT} and \eqref{equation:ergodicityexp} (see \cite{Hasminskii80sso} in the more general context of SDEs with multiplicative noise).
\begin{ass} \label{Assumption:growth}
The vector field $f:\R^d\rightarrow \R^d$ is globally Lipchitz and $C^\infty$, and there exist $C_1,C_2>0$ such that for all $x\in\R^d$,
$$
x^T f(x) \leq -C_1 x^T x + C_2.
$$
\end{ass}
The following stronger assumption yields the important special case of Brownian dynamics \eqref{equation:Langevin}.
\begin{ass}
\label{Assumption:f_gradient}
The vector field $f$ is a globally Lipschitz gradient, i.e.\ts there exists a $C^\infty$ potential $V:\R^d\rightarrow \R$ such that $f(x)=-\nabla V(x)$ is globally Lipschitz and
there exist $C_1>0$ and $C_2$ such that for all $x\in\R^d$, $V(x) \geq C_1 x^T x - C_2$.
\end{ass}
We recall that under Assumption \ref{Assumption:f_gradient}, the density
of the unique invariant measure is given by 
$\rho_\infty=Z\exp\left(-\frac{2V}{\sigma^2}\right)$ where $Z$ 
 is such that  $\int_{\R^d} \rho_\infty(x)dx=1$ but $Z$ is not numerically known in general.
We note that $\nabla \rho_\infty =\rho_\infty\frac{2}{\sigma^2}f$, 
equivalently $\nabla(\log \rho_\infty) = \frac{2}{\sigma^2}f$.

We denote $\CC^\infty_P(\R^d,\R)$ the vector space of $\CC^\infty$ functions such that all partial derivatives $\phi$ up to all orders have a polynomial growth of the form
$$
|\phi(x)| \leq C(1+|x|^s)
$$
for some constants $s$ and $C$ independent of $x$ (but depending on the order of differentiation).
For $\phi \in \CC^\infty_P(\R^d,\R)$ we define $u(x,t)=\E[\phi(X(t))|X(0)=x]$.
A classical tool for the study of \eqref{equation:sde} is the backward Kolmogorov equation \cite[Chap.\ts 2]{Milstein04snf}, which states that $u(x,t)$ solves the following deterministic parabolic PDE in $\R^d$,
\begin{equation}\label{equation:Kolmogorov}
\frac{\partial u}{\partial t} = \LL u, \qquad u(x,0)=\phi(x),
\qquad x\in \R^d,t>0,
\end{equation}
where the generator $\LL$ is defined as
\begin{equation} \label{equation:generator}
\LL\phi = f\cdot \nabla \phi + \frac{\sigma^2}2 \Delta \phi
\end{equation}
where $\Delta\phi=\sum_{i=1}^d \frac{\partial^2\phi}{\partial x_i^2}$
denotes the Laplace operator.
We recall that, under Assumption \ref{Assumption:growth} or \ref{Assumption:f_gradient}, the density of the invariant measure satisfies
$$
\LL^*\rho_\infty =0
$$
where $\LL^*\phi = -\mathrm{div}(f\phi) + \frac{\sigma^2}2 \Delta \phi$ is the $L^2$--adjoint of $\LL$.

We make the following natural assumptions on the numerical integrator \eqref{equation:defPsi}.

\begin{ass}
\label{Assumption:moments}
The numerical scheme \eqref{equation:defPsi} has bounded moments of any order along time, i.e.\ts for all integer $k\geq 0$,
$$
\sup_{n\geq 0} \E[|X_n|^{2k}] <\infty.
$$
\end{ass}

\begin{ass}
\label{Assumption:A0=L}
The numerical scheme \eqref{equation:defPsi} has a weak Taylor expansion of the form
\begin{equation}
\label{equation:dvp_U}
\E[\phi(X_1)|X_0=x] = \phi(x) + h\AA_0 \phi(x) + h^2 \AA_1\phi(x)+ \cdots
\end{equation}
for all $\phi \in \CC_P^\infty(\R^d,\R)$,
where $\AA_i,~i=0,1,2,\ldots$ are linear differential operators with coefficients depending smoothly on the drift $f$, and its derivatives (and depending on the choice of the integrator). In addition, we assume that $\AA_0$
coincides with the generator $\mathcal{L}$ given in \eqref{equation:generator},
i.e.\ts it is consistent and has (at least) local order one in the weak sense, $\AA_0 = \mathcal{L}$.
\end{ass}

\begin{remark}
A convenient sufficient condition to satisfy Assumption \ref{Assumption:moments} is given in \cite[Lemma 2.2.2]{Milstein04snf}: if $X_0$ is deterministic or has bounded moments of all order and the Markov chain $(X_n)_n$ satisfies
$$\abs{\E[X_{n+1}-X_n|X_n]}\leq C(1+\abs{X_n})h, \qquad \abs{X_{n+1}-X_n}\leq M_n(1+\abs{X_n})\sqrt{h},
$$
for $C$ a constant independent of $h$ and $M_n$ a random variable whose moments are all bounded uniformly with respect to $h$ small enough, then the numerical scheme satisfies Assumption \ref{Assumption:moments}.
Note that consistent Runge-Kutta type schemes such as \eqref{equation:defRK} satisfy Assumption \ref{Assumption:moments} (using the global Lipschitz assumption of $f$) and Assumption \ref{Assumption:A0=L}.
\end{remark}

In the following theorem, we recall the characterisation of order $p$ for the invariant measure of an ergodic integrator in terms of adjoints of the linear differential operators of Assumption \ref{Assumption:A0=L}. It was shown in \cite{Abdulle14hon} in the more general context of ergodic SDEs with multiplicative noise based on backward error analysis results in \cite{Debussche12wbe} on the torus
and generalizations \cite{Kopec15wbea,Kopec15wbeb} in the space $\R^d$.
\begin{theorem} \cite{Abdulle14hon}
\label{thm:orderAstar}
Assume Assumption \ref{Assumption:growth} or \ref{Assumption:f_gradient}. Consider the one step integrator \eqref{equation:defPsi} and assume that it is ergodic when applied to \eqref{equation:sde}.
Assume further Assumptions \ref{Assumption:moments} and \ref{Assumption:A0=L}.
If
\begin{equation} \label{equation:Ajstar}
\AA_j^* \rho_\infty = 0, \qquad j=2,\ldots ,p-1,
\end{equation}
then the scheme has order $p$ for the invariant measure and \eqref{equation:order}-\eqref{equation:orderexp} hold.
\end{theorem}
\begin{remark} 
Assuming in addition that the scheme has weak order $p-1$ of accuracy, i.e.\ts assuming in addition the stronger assumption 
$
\AA_j = \LL^{j+1}/(j+1)!, j=2,\ldots ,p-2
$,
then Theorem \ref{thm:orderAstar} is an immediate consequence of the
Talay-Tubaro expansion of the error \cite{Milstein85wao,Talay90eot} (see also \cite[Chap.\ts 2.2,\ts 2.3]{Milstein04snf}) given by
$$
\int_{\R^d} \phi(x) d\mu^h(x) - \int_{\R^d} \phi(x) d\mu(x) = \lambda_p h^p + \OO(h^{p+1})
$$
with 
$$
\lambda_{p}=\int_{0}^{+\infty}\int_{\R^{d}}\left({A}_{p}-\frac{1}{(p+1)!}\mathcal{L}^{p+1}\right)u(y,t) \rho_\infty(y)dydt
$$
where $u(x,t)$ is the solution of \eqref{equation:Kolmogorov}. 
Indeed, considering the $L^2$--adjoint of the operator ${A}_{p}-\frac{1}{(p+1)!}\mathcal{L}^{p+1}$ and using $\mathcal{L}^*\rho_\infty=0$ and $\AA_p^*\rho_\infty=0$ then yield $\lambda_p=0$ and the scheme has (at least) order $p$ for the invariant measure.
\end{remark}

The following extension of Theorem \ref{thm:orderAstar} permits 
to combine an integrator \eqref{equation:defPsi} with a postprocessor to achieve
high order for the invariant measure at a negligible overcost compared to a standard scheme.
\begin{theorem} 
\label{Theorem:Postprocessed_scheme_order_condition}
\cite{Vilmart15pif}
Assume the hypotheses of Theorem \ref{thm:orderAstar} and consider a postprocessor
$$
\overline X_n = G_n(X_n)
$$
that admits the following weak Taylor expansion for all $\phi\in C_P^\infty(\R^d,\R)$,
\begin{equation} \label{equation:dvp_postproc}
\E[\phi(G_n(x))]=\phi(x)+\sum_{i=1}^{p-1} \alpha_i h^i \LL^i \phi(x)+h^p\overline{\AA_p}\phi(x)+\cdots,
\end{equation}
for some constants $\alpha_i$ and a linear differential operator $\overline{\AA_p}$.
Assume further that
\begin{equation} \label{equation:Ajstarproc}
(\AA_p+[\LL,\overline \AA_p])^* \rho_\infty=0
\end{equation}
where $[\LL,\overline \AA_p]=\LL\overline \AA_p-\overline \AA_p\LL$ is the Lie bracket.
Then $\overline X_n$ yields an approximation of order $p+1$ for the invariant measure and it satisfies \eqref{equation:order} and \eqref{equation:orderexp}
with $p$ replaced by $p+1$ and $X_n$ replaced by $\overline X_n$.
\end{theorem}

Theorem \ref{Theorem:Postprocessed_scheme_order_condition}
 is stated and proved in \cite{Vilmart15pif} in the special case $\alpha_i=0$, $i=1,\ldots,p-1$ in \eqref{equation:dvp_postproc}.
However, the proof for non zero $\alpha_i$'s is nearly identical and thus is omitted.
Notice that the order conditions \eqref{equation:Ajstar} and \eqref{equation:Ajstarproc} are respectively equivalent to the identities
\begin{align*}
\int_{\R^d} (\AA_j \phi) \rho_\infty dx&=0,\qquad j=2,\ldots,p-1,\\
\int_{\R^d} (\AA_p \phi+[\LL,\overline \AA_p] \phi)  \rho_\infty dx&=0,
\end{align*}
for all test function $\phi\in C_P^\infty(\R^d,\R)$.
In the following section, we introduce the suitable algebraic framework based on exotic aromatic trees and B-series for the systematic study of these order conditions of accuracy for the invariant measure.

\section{Exotic aromatic trees and forests}
\label{section:exotic_aromatic_forests}

We first recall the known framework of aromatic B-series before introducing a modification well suited for invariant measure order conditions and called exotic aromatic B-series. We rely on the aromatic trees and forests introduced in \cite{Chartier07pfi} and rely on the presentation in \cite{Bogfjellmo15aso}.

\subsection{Aromatic trees and forests}
\label{section:aromatic_forests}

We first consider directed graphs $\gamma=(V,E)$ with $V$ a finite set of nodes and $E\subset V\times V$ the set of directed edges. If $(v,w)\in E$, we say that the edge is going from $v$ to $w$, and $v$ is called a predecessor of $w$.
Two directed graphs $(V_1,E_1)$ and $(V_2,E_2)$ are equivalent if there exists a bijection $\varphi: V_1\rightarrow V_2$ with $(\varphi\times\varphi) (E_1)=E_2$.
For brevity of notation, to avoid drawing arrows on the forests, an edge linking two nodes goes from the top node to the bottom one. If there is an eventual cycle, the arrows on it are going in the clockwise direction. For example, $$\includegraphics[scale=0.5]{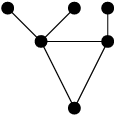}=\includegraphics[scale=0.5]{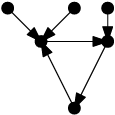}.$$

We call aromatic forests the equivalence classes of directed graphs where each node has at most one outgoing edge.
The connected components making an aromatic forest are called aromatic trees.
According to the above definition, there are two types of trees:
\begin{itemize}
\item \emph{aromas} are aromatic trees\footnote{Such graphs with one cycle are not strictly speaking ``trees'', they are however called aromatic trees in the literature as an analogy with carbon chemistry.} with exactly one cycle: $\etree 1 2 0 1$, $\etree 2 2 0 1$, $\etree 5 2 0 1$, $\etree 5 2 0 2$, \dots
\item \emph{rooted trees} do not have a cycle ; they have a unique node that has no outgoing edge and that is called the root, graphically represented at the bottom: $\etree 1 1 0 1$, $\etree 2 1 0 1$, $\etree 3 2 0 1$, $\etree 3 1 0 1$, \dots
\end{itemize}
Thus, an aromatic forest is a collection of aromas and rooted trees.
We call $\AA\TT=\{\etree 1 1 0 1, \etree 1 2 0 1 \: \etree 1 2 0 1 \: \etree 2 1 0 1, \dots\}$ the set of aromatic forests containing exactly one rooted tree, and we name its elements the aromatic rooted forests.

\begin{definition}[Elementary differentials]
\label{Elementary_differentials}
Let $\gamma=(V,E)\in \AA\TT$, and let $f: \R^d \rightarrow \R^d$ be a smooth function. We denote $\pi(v)=\{w\in V, (w,v)\in E\}$ the set of all predecessors of the node $v\in V$ and $r$ the root of $\gamma$. We also call $V^0=V\smallsetminus\{r\}=\{v_1,\dots ,v_m\}$ the other nodes of $\gamma$.
Finally we introduce the notation $I_{\pi(v)}=(i_{q_1},\dots ,i_{q_s})$ where the $q_k$ are the predecessors of $v$, and $$\partial_{I_{\pi(v)}} f=\frac{\partial^s f}{\partial x_{i_{q_1}}\dots \partial x_{i_{q_s}}}.$$
Then $F(\gamma)$ is defined as $$F(\gamma)(f)=\sum_{i_{v_1},\dots ,i_{v_m}=1}^d \left(\prod_{v\in V^0} \partial_{I_{\pi(v)}} f_{i_v}\right)\partial_{I_{\pi(r)}} f.$$
\end{definition}

\begin{ex*}
Let $\gamma=\includegraphics[scale=0.5]{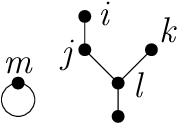}$ and $\widetilde{\gamma}=\includegraphics[scale=0.5]{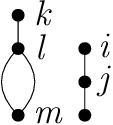}$ in $\AA\TT$
where we added indices to apply the formula of Definition \ref{Elementary_differentials}. Note that there is no index for the root.
Then the associated differentials are respectively
$F(\gamma)(f)=\sum_{i,j,k,l,m=1}^d \partial_{m} f_{m} f_i \partial_{i} f_{j} f_k \partial_{j,k}f_l \partial_l f=\Div(f)\cdot f'f''(f'f,f)$ and $F(\widetilde{\gamma})(f)=\sum_{i,j,k,l,m=1}^d \partial_{l} f_{m} \partial_{m,k} f_{l} f_k f_i \partial_{i} f_{j} \partial_{j}f=\sum_{m=1}^d f'_m((\partial_m f)'(f)) \cdot f'f'f$.
\end{ex*}

\subsection{Exotic aromatic trees and forests}
\label{section:eat}

We now introduce a new kind of edge, called a liana, for the aromatic forests. The corresponding generalization is called exotic aromatic forests.
Let $(V,E)$ be an aromatic forest and $L$ be a finite list of pairs of elements of $V$ (possibly with duplicates), then $\gamma=(V,E,L)$ is an exotic aromatic forest.
The elements of $L$ are called lianas and correspond to non-oriented edges between any two nodes of the forest. We graphically represent them with a dashed edge linking the two given nodes. As we authorize duplicates, there can be several lianas between two given nodes. Also lianas can link a node to itself.
For a node $v$, $\Gamma(v)$ denotes the list of the lianas (also with possible duplicates) linked to $v$. The predecessors of $v$ only take in account the edges of $E$.
An exotic aromatic tree of an exotic aromatic forest $\gamma=(V,E,L)$ is a connected component of the associated aromatic forest $(V,E)$.
We call $\EE\AA\TT$ the set of exotic aromatic forests with exactly one rooted tree, and name its elements exotic aromatic rooted forests.

\begin{ex*}
The lianas can link different trees of an aromatic forest and thus yield an exotic aromatic forest.
For instance, linking the aroma $\etree 1 2 0 1$ and the rooted tree $\etree 5 3 0 1$ gives $\includegraphics[scale=0.5]{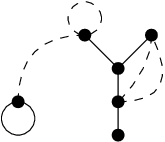}$.
\end{ex*}

The definition of elementary differentials is extended as follows.

\begin{definition}
\label{elementary_differential}
Let $\gamma=(V,E,L)\in \EE\AA\TT$, and let $f: \R^d \rightarrow \R^d$ be a smooth function. We name $r$ the root of $\gamma$ and $V^0=V\smallsetminus\{r\}=\{v_1,\dots ,v_m\}$ the other nodes of $\gamma$. We denote $l_1$,\dots ,$l_s$ the elements of $L$ and for $v\in V$, $J_{\Gamma(v)}$ the multiindex $(j_{l_{x_1}},\dots ,j_{l_{x_t}})$ where $\Gamma(v)=\{l_{x_1},\dots ,l_{x_t}\}$.
Then $F(\gamma)$ is defined as $$F(\gamma)(f)=\sum_{i_{v_1},\dots ,i_{v_m}=1}^d \sum_{j_{l_1},\dots ,j_{l_s}=1}^d \left(\prod_{v\in V^0} \partial_{I_{\pi(v)}} \partial_{J_{\Gamma(v)}} f_{i_v}\right) \partial_{I_{\pi(r)}} \partial_{J_{\Gamma(r)}} f.$$
\end{definition}

\begin{exs*}
The differential that corresponds to the rooted tree $\etree 1 1 1 1$ with a single node and a single liana is $F(\etree 1 1 1 1)(f)=\Delta f$.
We can also represent as exotic aromatic forest more complicated derivatives. For instance, let $\gamma=\includegraphics[scale=0.5]{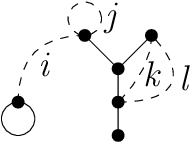}$,
then 
$$F(\gamma)(f)=\sum_{i,j,k=1}^d \Div(\partial_i f)\cdot f'((\partial_{kl} f)'(f''(\partial_{ijj} f,\partial_{kl} f))).$$
\end{exs*}

\subsection{Grafted exotic aromatic trees}

For the study of the order for the invariant measure of numerical integrators, we introduce an extension of exotic aromatic forests. The root now symbolizes a test function $\phi$, and it has leafs (nodes without predecessors) that represent a random standard normal vector $\xi$. Note that these new trees can be seen as bi-coloured trees in the context of P-series (see \cite[Chap.\ts 3]{Hairer06gni}), where the nodes represented with crosses cannot have predecessors.

\begin{definition}
\label{definition:grafted_trees}
A grafted node is a new type of node graphically represented by a cross.
Let $V$ be a set of nodes whose subset of grafted nodes is $V_g$, let $E$ be a set of edges such that each node in $V_g$ has exactly one outgoing edge and no ingoing edge, and let $L$ be a set of lianas that link nodes in $V\smallsetminus V_g$, then $\gamma=(V,E,L)$ is a grafted exotic aromatic forest.
We define as before the grafted exotic aromatic trees and grafted exotic aromatic rooted forests, that we denote $\EE\AA\TT_g$. 

If $\gamma=(V,E,L)$ is a grafted exotic aromatic rooted forest, $\phi: \R^d\rightarrow \R$ a smooth function, and $\xi$ a random vector of $\R^d$ whose components are independent and follow a standard normal law, the associated elementary differential of $\gamma$ is, with the same notation as Definition \ref{elementary_differential} and $V^0=V\smallsetminus (V_g\cup \{r\})$,
$$F(\gamma)(f,\phi,\xi)=\sum_{i_{v_1},\dots ,i_{v_m}=1}^d \sum_{j_{l_1},\dots ,j_{l_s}=1}^d \left(\prod_{v\in V^0} \partial_{I_{\pi(v)}} \partial_{J_{\Gamma(v)}} f_{i_v}\right)
\left(\prod_{v\in V_g} \xi_{i_v}\right) \partial_{I_{\pi(r)}} \partial_{J_{\Gamma(r)}} \phi.$$
\end{definition}

\begin{ex*}
The differential associated to the forest $\ctree 2 1 2 1$ is $F(\ctree 2 1 2 1)(f,\phi,\xi)=\phi'(f''(\xi,\xi))$.
\end{ex*}

If $\gamma$ is such that $V_g$ is empty, we recover the exotic aromatic forests of Definition \ref{elementary_differential}, where $\phi$ is replaced by $f$. For the rest of the paper (except Section \ref{section:isometric_equivariance}), we update the definition of the elementary differential of an exotic aromatic forests so that the root is associated to the function $\phi$.
This definition can be straightforwardly extended on non-rooted exotic aromatic forests.
For brevity of notation, we also write $F(\gamma)(\phi)$ instead of $F(\gamma)(f,\phi,\xi)$. We note that $\phi \to F(\gamma)(\phi)$ is a linear differential operator (dependent on $f$ and $\xi$).

\subsection{Grafted exotic aromatic B-series}

In this section, we adapt the formalism of aromatic B-series of \cite{MuntheKaas16abs} to grafted exotic aromatic forests, in order to use it as a numerical tool for weak Taylor expansions in the next sections.
We define the order $\abs{\gamma}$ of a tree $\gamma \in \EE\AA\TT_g$.
We denote $N(\gamma)$ the number of nodes, $N_l(\gamma)$ the number of lianas, $N_c(\gamma)$ the number of grafted nodes and $N_v(\gamma)=N(\gamma)-N_c(\gamma)-1$ the number of nodes that are non grafted and different from the root, then
$$\abs{\gamma}=N_v(\gamma)+N_l(\gamma)+\frac{N_c(\gamma)}{2}.$$

\begin{definition}
Let $a: \EE\AA\TT_g\to \R$ a map, and let $f: \R^d \rightarrow \R^d$ and $\phi: \R^d \rightarrow \R$ be two smooth functions, then the grafted exotic aromatic B-series $B(a)(\phi)$ is a formal series indexed over $\EE\AA\TT_g$ defined by
$$B(a)(\phi)=\sum_{\gamma \in \EE\AA\TT_g} h^{\abs{\gamma}} a(\gamma) F(\gamma)(\phi).$$
We extend the definition of $F$ on $\Vect(\EE\AA\TT_g)$ by writing
$$F\lpar\sum_{\gamma \in \EE\AA\TT_g} h^{\abs{\gamma}}a(\gamma) \gamma\rpar(\phi)=B(a)(\phi).$$
\end{definition}

The variable $h$ is formal and thus can be chosen to be equal to 1.
If the series is indexed only on (exotic) aromatic rooted forests, then it is called an (exotic) aromatic B-series. In Section \ref{section:isometric_equivariance}, we shall focus on exotic aromatic B-series.

\begin{remark}
The coefficients $a(\gamma)$ of standard B-series are sometimes renormalized as $\frac{a(\gamma)}{\rho(\gamma)}$ where $\rho$ is a function determined by the symmetries of the associated forest. If $\rho$ is appropriately chosen, it greatly simplifies the composition laws of (aromatic) B-series (see \cite{Hairer06gni,Chartier10aso,Bogfjellmo15aso}).
Finding the best definition of $\rho$ for this exotic extension of B-series is out of the scope of this paper.
\end{remark}

\subsection{Isometric equivariance of exotic aromatic rooted forests}
\label{section:isometric_equivariance}

In this subsection, we show that the exotic aromatic B-series satisfy an isometric equivariance property in the spirit of \cite{MuntheKaas16abs, McLachlan16bsm}.
We consider exotic aromatic rooted forests $\gamma$ where the differential associated to the root is $f$. As the function $f$ is no longer fixed, we denote the associated differential $F(\gamma)(f)$. Also we adapt the definition of exotic aromatic B-series to this change.
First we add a new tree: the empty tree $\varnothing$.
The function $F$ is then extended on $\EE\AA\TT_g\cup\{\varnothing\}$ by $F(\varnothing)(f)=\Id_{\R^d}$.
Then, for a function $a:\EE\AA\TT\cup\{\varnothing\}\rightarrow \R$, the associated exotic aromatic B-series is
$$B(a)(f)=\sum_{\gamma \in \EE\AA\TT\cup\{\varnothing\}} a(\gamma) F(\gamma)(f).$$
We study (exotic) aromatic B-series $B(a)$ with $a(\varnothing)=1$. We call these (exotic) aromatic B-series methods.
Let $G$ be a subgroup of $\GL_d(\R)\ltimes \R^d$, let the action of an element $(A,b)\in G$ on $\R^d$ be
$x\mapsto Ax+b$,
and let the action on a vector field $f:\R^d \to \R^d$ be
$$((A,b)*f)(x):=Af(A^{-1}(x-b)).$$
We simplify the notation by writing $A*f:=(A,0)*f$.
We recall the definition of equivariance from \cite{MuntheKaas16abs}. The property of equivariance means the method is unchanged by an affine coordinate transformation.
Let $\Phi$ be a differential operator and let $G$ be a subgroup of $\GL_d(\R)\ltimes \R^d$, then $\Phi$ is called $G$-equivariant if
$$\forall (A,b)\in G, \forall f\in \CC^\infty(\R^d,\R^d), \Phi((A,b)*f)=(A,b)\circ\Phi(f)\circ (A,b)^{-1}.$$
In particular, $\Phi$ is said to be affine equivariant if $G=\GL_d(\R)\ltimes \R^d$ and isometric equivariant if $G=\Or_d(\R)\ltimes \R^d$.
\begin{theorem}
\label{theorem:equivariance}
Consider an exotic aromatic B-series method $B(a)$ (with $a(\varnothing)=1$), then $B(a)$ is isometric equivariant.
\end{theorem}
\begin{remark}
It is proved in \cite{MuntheKaas16abs} that standard B-series methods are
exactly the affine equivariant methods.
Analogously, it would be interesting to characterize the isometric equivariant maps.
\end{remark}
For the sake of brevity, we omit the proof of Theorem \ref{theorem:equivariance}.
The proof can be made in 
the spirit of the result \cite[Prop.\ts 2.1]{MuntheKaas16abs} for affine equivariant B-series.

\section{Analysis of invariant measure order conditions using exotic aromatic forests}
\label{section:using_eat}

In this section, we show how the framework of Section \ref{section:exotic_aromatic_forests} applies for the study of order conditions for the invariant measure of numerical integrators.

\subsection{Weak Taylor expansion using exotic aromatic forests}
Let us begin this section with the example of the $\theta$-method \eqref{equation:theta_method}.
We apply the usual methodology to expand in Taylor series $\E[\phi(X_1)|X_0=x]$ as $h\to 0$.
We refer to \cite{Zygalakis11ote,Abdulle12hwo} for other examples of analogous calculations performed without exotic aromatic forests.
Under $X_0=x$, we have $$X_1=x+\sqrt{h}\sigma \xi+hf+h\sqrt{h}\theta \sigma f'\xi+h^2\theta f'f+h^2 \frac{\theta\sigma^2}{2} f''(\xi,\xi)+\cdots  $$
Then we deduce $\E[\phi(X_1)|X_0=x]=\phi(x)+ h\LL\phi (x) +h^2 \AA_1 \phi (x) +\cdots $,
where 
\begin{align}
\label{equation:A1_with_crosses}
\AA_1\phi & = \E[\theta \phi'f'f+\frac{1}{2}\phi''(f,f)+\frac{\theta \sigma^2}{2}\phi'f''(\xi,\xi) +\theta \sigma^2 \phi''(f'\xi,\xi) \nonumber\\
 & +\frac{\sigma^2}{2} \phi^{(3)}(f,\xi,\xi) +\frac{\sigma^4}{24} \phi^{(4)}(\xi,\xi,\xi,\xi)] \nonumber\\
 & = \mathbb{E}[F(\theta \etree 3 1 0 1+\frac{1}{2}\etree 3 2 0 1+\frac{\theta \sigma^2}{2}\ctree 2 1 2 1 +\theta \sigma^2 \ctree 2 1 2 2 +\frac{\sigma^2}{2} \ctree 2 1 2 3 +\frac{\sigma^4}{24} \ctree 1 1 4 1 )(\phi)].
\end{align}
All the forests with an odd number of grafted nodes vanished because odd moments of a centred Gaussian random variable are zero.
The expectation of the differential of a forest with exactly two grafted nodes comes straightforwardly.
$$
\mathbb{E}[F(\ctree 2 1 2 1 )(\phi)]=\E[\phi'(f''(\xi,\xi))]
=\sum_{i,j,k} \partial_i \phi  \partial_{jk} f_i  \E[\xi_j\xi_k]
=\sum_{i,j} \partial_i \phi  \partial_{jj} f_i
=F(\etree 2 1 1 1)(\phi),
$$
where $\E[\xi_j\xi_k]=0$ for $j\neq k$ because of the independence of the $\xi_i$'s.
We see that taking the expectation of the differential associated to a grafted tree amounts to linking the grafted nodes with lianas in all possible manners.
For instance, for the following example with four grafted nodes:
\begin{align*}
\mathbb{E}[F(\ctree 1 1 4 1 )(\phi)]&=\E[\phi^{(4)}(\xi,\xi,\xi,\xi)]
=\sum_{i,j,k,l} \partial_{i,j,k,l} \phi  \E[\xi_i\xi_j\xi_k\xi_l]\\
&=\sum_{i} \partial_{i,i,i,i} \phi  \E[\xi_i^4]+3\sum_{\underset{i\neq j}{i,j}} \partial_{i,i,j,j} \phi  \E[\xi_i^2]\E[\xi_j^2]
=3\sum_{i,j} \partial_{i,i,j,j} \phi\\
&=3F(\etree 1 1 2 1)(\phi).
\end{align*}
Let us now comment this computation. The interesting fact is that $\E[\xi_i^4]=3$ corresponds exactly to the number of ways to gather the indices $i$, $j$, $k$ and $l$ in pairs. This observation makes an exotic aromatic tree naturally appear.
However, here we took only four grafted nodes and the differential form was symmetric in the arguments $\xi$. We need to study the expectation of general exotic aromatic forest elementary differentials.
This is the aim of the following theorem.

\begin{theorem}
\label{Theorem:Forest_expectation}
Let $\gamma\in \EE\AA\TT_g$ be a grafted exotic aromatic rooted forest with an even number of grafted nodes $2n$, let $\phi:\R^d\rightarrow \R$ be a smooth function, and let $V^\times=\{c_1,\dots ,c_{2n}\}$ be the set of grafted nodes of $\gamma$. We call $\PP_2(2n)$ the set of partitions by pair of $\{1,\dots ,2n\}$, i.e.\ts the set of surjections $p: \{1,\dots ,2n\}\rightarrow \{1,\dots ,n\}$ such that the preimage of each singleton has exactly two elements and the minima of those preimages follow an ascending order ($\min(p^{-1}(\{i\}))<\min(p^{-1}(\{j\}))$ for $i<j$).
Finally we define $\varphi_\gamma: \PP_2(2n) \rightarrow \EE\AA\TT$ the application that maps the partition $p$ of $\gamma$ to the aromatic forest where the grafted nodes are linked by lianas according to $p$.
Then, the expectation of $F(\gamma)(\phi)$ is given by
$$
\E[F(\gamma)(\phi)]=\sum_{p\in\PP_2(2n)} F(\varphi_\gamma(p))(\phi).
$$
\end{theorem}

This theorem states that taking the expectation of the differential associated to a forest amounts to sum the forests obtained by linking the grafted nodes together pairwise using lianas in all possible manners and take the associated differential.

\begin{ex*}
Let us take $\gamma \in \{\etree 1 1 0 1, \ctree 1 1 2 1, \ctree 1 1 4 1, \dots \}$ the tree with only a root, $2n$ grafted nodes and no liana, then
$$\E[F(\gamma)(\phi)]=\frac{(2n)!}{2^nn!}\Delta^n \phi$$
The integer $\frac{(2n)!}{2^nn!}$ is exactly the number of ways to gather the grafted nodes by pairs.
An other example is
$$\E[F(\ctree 2 1 4 1)(\phi)]=3 F(\etree 2 1 2 4),$$
where the coefficient 3 accounts for the number of choices for linking the grafted nodes pairwise.
\end{ex*}

\begin{app*}
Using Theorem \ref{Theorem:Forest_expectation}, we immediately obtain that Runge-Kutta methods \eqref{equation:defRK} can be developed in exotic aromatic forests.
As a special Runge-Kutta method, we get back to the $\theta$-method \eqref{equation:theta_method}.
The operator $\AA_1$ is now convenient to write with exotic aromatic trees. Applying Theorem \ref{Theorem:Forest_expectation} to \eqref{equation:A1_with_crosses}, we deduce $\AA_1=F(\gamma)$ with
\begin{equation}
\label{equation:gamma_theta_method}
\gamma=\theta \etree 3 1 0 1
+\frac{1}{2}\etree 3 2 0 1
+\frac{\theta \sigma^2}{2}\etree 2 1 1 1
+\theta \sigma^2 \etree 2 1 1 2
+\frac{\sigma^2}{2} \etree 2 1 1 3
+\frac{\sigma^4}{8} \etree 1 1 2 1 .
\end{equation}
\end{app*}

Theorem \ref{Theorem:Forest_expectation} follows from the following lemma, which is an extension of the Isserlis theorem \cite{Isserlis18oaf} to the case of multilinear mappings.
The Isserlis theorem states that if $\chi$ is a $2n$-dimensional Gaussian random vector with mean zero and arbitrary covariance, then $$\E\left[\prod_{i=1}^{2n} \chi_i\right]=\sum_{p\in\PP_2(2n)} \prod_{\underset{p(i)=p(j)}{i<j}} \E[\chi_i\chi_j].$$
For $n=2$, it gives $\E[\chi_1 \chi_2 \chi_3 \chi_4]=
\E[\chi_1 \chi_2]\E[\chi_3 \chi_4]
+\E[\chi_1 \chi_3]\E[\chi_2 \chi_4]
+\E[\chi_1 \chi_4]\E[\chi_2 \chi_3]$.

\begin{lemma}
\label{Lemma:Forest_expectation}
Let $B:\R^d\times\dots \times \R^d=\R^{2nd}\to \R$ be a $2n$-multilinear form, and let $\xi$ be a Gaussian vector $\NN(0,I_d)$, then
\begin{equation}
\label{equation:lemme_expectation}
\E[B(\xi,\dots ,\xi)]=\sum_{p\in\PP_2(2n)} \sum_{i_{1},\dots ,i_n=1}^d B(e_{i_p}),
\end{equation}
with $e_{i_p}=(e_{i_{p(1)}},\dots , e_{i_{p(2n)}})$, and we recall that $e_1,\ldots,e_{d}$ denotes the canonical basis of $\R^d$.
\end{lemma}

\begin{proof}
For the particular case of an elementary multilinear form $B_{\sigma}:(x_1,\dots ,x_{2n})\mapsto \prod_{j=1}^{2n} (x_j)_{\sigma(j)}$ where $\sigma : \{1,\dots ,2n\}\rightarrow \{1,\dots ,d\}$ is a given mapping and $(x_j)_{\sigma(j)}$ denotes the $\sigma(j)$'s component of $x_j\in \R^d$, the identity \eqref{equation:lemme_expectation} reduces to the Isserlis theorem. As any multilinear form can be decomposed as a linear combination of such elementary multilinear forms, the result \eqref{equation:lemme_expectation} is proved by linearity with respect to $B$.
\end{proof}

\begin{proof}[Proof of Theorem \ref{Theorem:Forest_expectation}]
We consider $F(\gamma)(\phi)$ as a $2n$-multilinear form $B_{\gamma,\phi}$ evaluated in $(\xi,\dots ,\xi)$ (see Definition \ref{definition:grafted_trees}). Lemma \ref{Lemma:Forest_expectation} gives
$$\E[B_{\gamma,\phi}(\xi,\dots ,\xi)]=\sum_{p\in\PP_2(2n)} \sum_{i_{1},\dots ,i_n=1}^d B_{\gamma,\phi}(e_{i_p})=\sum_{p\in\PP_2(2n)} F(\varphi_\gamma(p))(\phi),$$
because $B_{\gamma,\phi}(e_{i_p})$ is the differential $F(\gamma)$ where we differentiate the non-grafted nodes linked to grafted nodes in the directions given by $e_{i_p}$, and $F(\varphi_\gamma(p))(\phi)$ is obtained by summing $B_{\gamma,\phi}(e_{i_p})$ over all the indices.
\end{proof}

\subsection{Integration by parts of the exotic aromatic forests}
\label{Section:IPP}
The goal of this section is to integrate by parts $\int_{\R^d} F(\gamma)(\phi)\rho_\infty dx$, for $\gamma$ an exotic aromatic rooted forest, in order to write it in the form $\int_{\R^d} \phi'\widetilde{f}\rho_\infty dx$ for a certain sum of elementary differentials $\widetilde{f}$. The idea is to transform a high order differential operator $\AA : \phi \to F(\gamma)(\phi)$ into a differential operator $\phi \to \phi'\widetilde{f}$ of order 1 such that $\AA^*\rho_\infty=-\Div(\widetilde{f}\rho_\infty)$. The tree formalism previously defined makes this task systematic and very convenient.
This also serves as a crucial ingredient in the next section.
Let us first begin with an example.
\begin{equation}
\label{equation:exemple_IPP}
\begin{split}
\int_{\R^d} F(\etree 3 2 0 1)(\phi) \rho_\infty  dx&=\sum \limits_{i,j} \int_{\R^d} \frac{\partial ^2 \phi}{\partial x_i \partial x_j}f_i f_j \rho_\infty  dx\\
&=-\sum \limits_{i,j} \lcro
\int_{\R^d} \frac{\partial \phi}{\partial x_j} \frac{\partial f_i}{\partial x_i} f_j \rho_\infty dx 
+ \int_{\R^d} \frac{\partial \phi}{\partial x_j} f_i \frac{\partial f_j}{\partial x_i} \rho_\infty  dx\\
&+ \int_{\R^d} \frac{\partial \phi}{\partial x_j} f_i f_j \frac{\partial \rho_\infty}{\partial x_i} dx
\rcro,
\end{split}
\end{equation}
where we integrated by parts; note that the boundary term vanishes using the growth assumptions on $\phi$.

\begin{notation*}
We denote $g=\log(\rho_\infty)$, then $\nabla \rho_\infty=(\nabla g) \rho_\infty$.
\end{notation*}

We have
$$\int_{\R^d} F(\etree 3 2 0 1)(\phi) \rho_\infty dx= -\int_{\R^d} \Div(f) \phi 'f \rho_\infty dx - \int_{\R^d} \phi'f'f \rho_\infty dx - \int_{\R^d} g'f \phi'f \rho_\infty dx.$$
By writing $\widetilde{f}=-(\Div(f)f+f'f+g'f f)$, we deduce 
$$\int_{\R^d} F(\etree 3 2 0 1)(\phi) \rho_\infty dx= \int_{\R^d} \phi'\widetilde{f}\rho_\infty dx.$$
We see that even for a simple forest, the integration by parts requires some calculations and a new term appears: the function $g=\log(\rho_\infty)$, and its derivatives. We use below the exotic aromatic forests to make this task easier.

\begin{definition}[Aromatic root and elementary differential]
An aromatic root is a new type of node represented by a square that has no outgoing edge. An exotic aromatic tree that has an aromatic root is considered as an aroma. The definition of the sets $\EE\AA\TT$ is extended to include these new aromas.

If $\gamma=(V,E,L)\in \EE\AA\TT$ is an exotic aromatic rooted forest whose set of aromatic roots is $V^1\subset V$, if $V^0=V\smallsetminus(V^1\cup\{r\})$ and $V=\{r,v_1,\dots,v_m\}$ where $r$ is the root of $\gamma$, if $l_1$,\dots,$l_s$ are the elements of $L$, then $F(\gamma)$ is defined as $$F(\gamma)(\phi)=\sum_{i_{v_1},\dots,i_{v_m}=1}^d \sum_{j_{l_1},\dots,j_{l_s}=1}^d
\left(\prod_{v\in V^1} \partial_{I_{\pi(v)}} \partial_{J_{\Gamma(v)}} g\right)
\left(\prod_{v\in V^0} \partial_{I_{\pi(v)}} \partial_{J_{\Gamma(v)}} f_{i_v}\right)
\partial_{I_{\pi(r)}} \partial_{J_{\Gamma(r)}} \phi.$$
\end{definition}

One can also extend the definition of $\EE\AA\TT_g$ in the same way so that it includes aromatic roots.

\begin{exs*}
$F(\gtree 2 1 0 1 \: \etree 2 1 0 1)(\phi)=g'f  \phi'f$ and $F(\gtree 2 2 1 1 \: \gtree 2 4 1 1)(\phi)=\sum_{i,j} (\partial_j g)^2 \partial_i g \partial_i \phi$.
\end{exs*}

Then the equality \eqref{equation:exemple_IPP} can be rewritten using forests:
$$\int_{\R^d} F(\etree 3 2 0 1)(\phi) \rho_\infty dx= -\int_{\R^d}  F(\etree 1 2 0 1 \: \etree 2 1 0 1)(\phi) \rho_\infty dx - \int_{\R^d} F(\etree 3 1 0 1)(\phi) \rho_\infty dx - \int_{\R^d} F(\gtree 2 1 0 1 \: \etree 2 1 0 1)(\phi) \rho_\infty dx.$$
We notice that integrating by parts $F(\etree 3 2 0 1)(\phi)$ amounts to unplugging an edge from the root and to replug it either to all the other nodes of the forests or to an aromatic root and then to sum over all possibilities.
This intuition is made rigorous in the following theorem.

\begin{theorem}
\label{Theorem:IPP}
Let $\gamma\in \EE\AA\TT$, we choose a direction to integrate by parts, i.e.\ts an edge or a liana $e$ connected to the root $r$.
Then $$\int_{\R^d} F(\gamma)(\phi) \rho_\infty dx=-\sum_{\widetilde{\gamma}\in U(\gamma,e)}\int_{\R^d} F(\widetilde{\gamma})(\phi) \rho_\infty dx,$$
where $U(\gamma,e)$ is the set of exotic aromatic rooted forests obtained by unplugging the chosen edge/liana $e$ and by linking it either to a node different from $r$ or to a new aromatic root.
\end{theorem}

\begin{remark}
Theorem \ref{Theorem:IPP} can be extended in the following way: if $n$ is a node of $\gamma$ and $e$ an edge ingoing to $n$ or a liana connected to $n$, then the same result holds if we replace $U(\gamma,e)$ by $U(\gamma,n,e)$, the set of exotic aromatic rooted forests obtained by unplugging the chosen edge/liana $e$ and by linking it either to a node different from $n$ or to a new aromatic root.
\end{remark}

\begin{proof}[Proof of Theorem \ref{Theorem:IPP}]
We call $r$ the root and $V^0=\{v_1,\dots,v_m\}$ the other nodes of $\gamma$. We suppose for simplicity that $\gamma$ does not have any aromatic root, otherwise the proof can be adapted straightforwardly. We denote $l_1$,\dots,$l_s$ the elements of $L$ and $v_{k_1}$,\dots,$v_{k_p}$ the elements of $\pi(r)$. We choose to integrate by parts in the direction of the edge $x_{i_{v_{k_1}}}$.

\begin{align*}
\int_{\R^d} F(\gamma)(\phi)\rho_\infty dx &=\sum_{i_{v_1},\dots,i_{v_m}=1}^d \sum_{j_{l_1},\dots,j_{l_s}=1}^d \int_{\R^d} \left(\prod_{v\in V^0} \partial_{I_{\pi(v)}} \partial_{J_{\Gamma(v)}} f_{i_v}\right) \partial_{i_{v_{k_1}}\dots i_{v_{k_p}}} \partial_{J_{\Gamma(r)}} \phi \rho_\infty dx\\
&=- \sum_{i_{v_1},\dots ,i_{v_m}=1}^d \sum_{j_{l_1},\dots ,j_{l_s}=1}^d \lcro \\
&\sum_{u\in V^0} \int_{\R^d} \left(\prod_{v\in V^0\smallsetminus \{u\}} \partial_{I_{\pi(v)}} \partial_{J_{\Gamma(v)}} f_{i_v}\right) \partial_{i_{v_{k_1}}} \partial_{I_{\pi(u)}} \partial_{J_{\Gamma(u)}} f_{i_u} \partial_{i_{v_{k_2}}\dots i_{v_{k_p}}} \partial_{J_{\Gamma(r)}} \phi \rho_\infty dx \\
& + \int_{\R^d} \left(\prod_{v\in V^0} \partial_{I_{\pi(v)}} \partial_{J_{\Gamma(v)}} f_{i_v}\right) \partial_{i_{v_{k_2}}\dots i_{v_{k_p}}} \partial_{J_{\Gamma(r)}} \phi \partial_{i_{v_{k_1}}} g \rho_\infty dx \rcro
\end{align*}
Each term of the sum on $u\in V^0$ is the differential associated to the forest $\gamma_u$. This forest is obtained by unplugging the root of its edge linking it to $v_{k_1}$, and sticking it to $u$.
The last term of the computation is the differential of the forest obtained by linking the unplugged edge to an aromatic root.
These terms are exactly what we expected, thus the theorem is proved for the case of edges.
For the case of integrating in the direction of a liana, the proof is nearly identical. We just need to develop $J_{\Gamma(r)}$ instead of $I_{\pi(r)}$.
\end{proof}

\begin{definition}
\label{definition:IPP_equivalence_relation}
Let $\gamma_1$ and $\gamma_2$ be two exotic aromatic B-series, we define the equivalence relation $\sim$ and write $\gamma_1 \sim \gamma_2$ if we can transform $\gamma_1$ into $\gamma_2$ by integrating by parts the associated differentials according to the procedure presented in the previous theorem.
\end{definition}

\begin{exs*}
The integration by parts (\ref{equation:exemple_IPP}) can now be simply rewritten as 
$$\etree 3 2 0 1 \sim -\etree 3 1 0 1-\gtree 2 1 0 1 \: \etree 2 1 0 1-\etree 1 2 0 1 \: \etree 2 1 0 1.$$
One can also iterate the process of integration by parts to fully simplify the trees. For example, we have
$$
\etree 2 1 1 2 \sim -\etree 2 1 1 1 -\gtree 3 3 1 1 \qquad
\gtree 3 3 1 2 \sim -\gtree 3 3 1 1 - \gtree 1 1 1 1 \: \etree 2 1 0 1 -\gtree 2 2 1 1 \: \etree 2 1 0 1
$$
and summing up yields
$$\etree 2 1 1 3 \sim -\etree 2 1 1 2 -\gtree 3 3 1 2 \sim \etree 2 1 1 1 +2\gtree 3 3 1 1 + \gtree 1 1 1 1 \: \etree 2 1 0 1 +\gtree 2 2 1 1 \: \etree 2 1 0 1.$$
Finally here is a last example that will be used in Section \ref{Section:High_Order_Integrators}. We apply the procedure of integration by parts to the forest $\etree 1 1 2 1$.
$$\etree 1 1 2 1 
\sim -\gtree 1 1 2 1 
\sim \gtree 1 1 1 1 \: \etree 1 1 1 1 + \gtree 2 2 1 1 \: \etree 1 1 1 1
\sim -\gtree 2 3 2 1 - \gtree 1 1 1 1 \: \gtree 2 4 1 1 - 2 \gtree 3 2 2 1 - \gtree 2 2 1 1 \: \gtree 2 4 1 1$$
Using analytic formulas, $\int_{\R^d} F(\etree 1 1 2 1)(\phi)\rho_\infty dx=\int_{\R^d} \phi'\widetilde{f}\rho_\infty dx$, where
$$\widetilde{f}_i=-\partial_i(\Delta g)-\Delta g \partial_i g-2\sum_{j=1}^d \partial_j g \partial_{i,j} g-\sum_{j=1}^d (\partial_j g)^2 \partial_i g.$$
\end{exs*}

\subsection{Order conditions using exotic aromatic forests}
\label{section:formal_ipp}
In this section, we adapt Theorem \ref{thm:orderAstar} in the context of exotic aromatic forests. In the spirit of traditional B-series, we give, under Assumption \ref{Assumption:f_gradient}, the general simplification for orders up to three of a general numerical method expandable in exotic aromatic B-series. With these, one can improve a method order as presented in Section \ref{Section:High_Order_Integrators}, or derive conditions on the method to achieve high order for the invariant measure as we do in Sections \ref{section:Runge_Kutta}, \ref{postprocessed integrators}, \ref{section:partitioned_methods} and \ref{Section:non_reversible_perturbation}.
It is worth noting that with only integration by parts, we can formally derive numerical methods (see Section \ref{Section:High_Order_Integrators} for an example), but under Assumption \ref{Assumption:f_gradient}, the methods can be simplified.

\begin{proposition}[Simplification rules]
\label{Proposition:Simplification_rule}
Under Assumption \ref{Assumption:f_gradient}, the two forest patterns gathered in each of the following pairs represent the same differential:
$$
\includegraphics[scale=0.5]{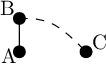} \text{ and } \includegraphics[scale=0.5]{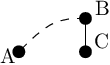}, \quad
\includegraphics[scale=0.5]{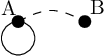} \text{ and } \includegraphics[scale=0.5]{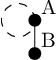}, \quad
\includegraphics[scale=0.5]{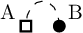} \text{ and } \frac{2}{\sigma^2} \includegraphics[scale=0.5]{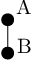}, \quad
\gtree 1 1 1 1 \text{ and } \frac{2}{\sigma^2} \etree 1 2 0 1.
$$
In the first and second cases, one can replace the nodes A, B, C with aromatic roots and the result remains. For the third case, the node B can also be replaced.
\end{proposition}

\begin{proof}
For the first pair of patterns, the associated differentials have the respective forms $\sum_{i=1}^d \partial_j f_k  \partial_i f_j  \partial_i f_l$ and $\sum_{j=1}^d \partial_j f_k  \partial_j f_i  \partial_i f_l$.
As $f$ is a gradient, $f'$ is a symmetric matrix and $\partial_i f_j=\partial_j f_i$. The two differentials are then equal.\\
The second point is proved in the same way. For the third and fourth points, we just use that $\nabla g=\frac{2}{\sigma^2}f$.
\end{proof}

\vspace{-6ex}

\begin{ex*}
We have $F(\gtree 2 3 2 1)=\frac{2}{\sigma^2}F(\etree 2 1 1 1)$ and $F(\gtree 3 3 1 1)=\frac{2}{\sigma^2}F(\etree 3 1 0 1)$.
\end{ex*}

The equivalence relation $\sim$ of Definition \ref{definition:IPP_equivalence_relation} is extended to include the simplification rules of Proposition \ref{Proposition:Simplification_rule}.

\begin{remark}
If $f$ is a general vector field not assumed to be a gradient, we can prove that the elementary differentials of exotic aromatic forests are independent.
The proof is an extension of the result in \cite[Chap.\ts 3, Exercise 3]{Hairer06gni}.
Indeed, for a given exotic aromatic forest $\gamma=(V,E,L)$, we take a bijective numbering $n_V:V\smallsetminus \{r\}\to \{1,2,\dots, \abs{V}-1\}$ of the nodes and another one for the lianas $n_L:L\to \{\abs{V},\dots, \abs{\gamma}\}$, then we define
\begin{align*}
\phi(x)&=\prod_{v\in \pi(r)} x_{n_V(v)} \prod_{l\in \Gamma(r)} x_{n_L(l)},\\
f_i(x)&=\prod_{v\in \pi(n_v^{-1}(i))} x_{n_V(v)} \prod_{l\in \Gamma(n_v^{-1}(i))} x_{n_L(l)}, \qquad i=1,\dots, \abs{V}-1,
\end{align*}
and $f_i=0$, $i=\abs{V},\dots, \abs{\gamma}$.
With this choice, $F(\widetilde{\gamma})(f, \phi)(0)\neq 0$ if and only if  $\widetilde{\gamma}=\gamma$, thus giving the independence of elementary differentials.
Note that this result does not hold if $f$ is assumed gradient due to
simplification rules (see Proposition \ref{Proposition:Simplification_rule}).
\end{remark}

We now adapt Theorem \ref{thm:orderAstar} to the context of exotic aromatic forests.

\begin{theorem}
\label{Theorem:Formal_IPP}
We assume Assumptions \ref{Assumption:growth} and \ref{Assumption:moments}. We consider an ergodic numerical scheme that can be developed in exotic aromatic B-series (and thus has a development of the form \eqref{equation:dvp_U})
\begin{equation}
\label{equation:dvp_U_Bseries}
\E[\phi(X_1)|X_0=x]=F(\etree 1 1 0 1)(\phi) + \sum_{\underset{1\leq\abs{\gamma}\leq p}{\gamma \in \EE\AA\TT}} h^{\abs{\gamma}} a(\gamma) F(\gamma)(\phi) + \cdots,
\end{equation}
for $p\in \N$. We denote $\AA_i=F(\gamma_i)$. If $\gamma_i\sim \widetilde{\gamma_i}$ and $F(\widetilde{\gamma_i})=0$ for all $1\leq i<p$, then the method is of order (at least) $p$ for the invariant measure. In particular, under Assumption \ref{Assumption:f_gradient}, we obtain for order 1
$$\widetilde{\gamma_0} = \lpar a\left(\etree 2 1 0 1\right)-\frac{2}{\sigma^2} a\left(\etree 1 1 1 1\right)\rpar \etree 2 1 0 1,$$
and in addition we have for order 2,
\begin{align*}
\widetilde{\gamma_1} &= \lpar  a(\etree 3 1 0 1)-\frac{2}{\sigma^2}a(\etree 2 1 1 2)+\frac{2}{\sigma^2}a(\etree 2 1 1 3)-\frac{4}{\sigma^4}a(\etree 1 1 2 1)\rpar \etree 3 1 0 1
+\lpar  a(\etree 2 1 1 1)-a(\etree 2 1 1 2)+a(\etree 2 1 1 3)-\frac{2}{\sigma^2}a(\etree 1 1 2 1)\rpar \etree 2 1 1 1\\
&+\lpar  a(\etree 3 2 0 1)-\frac{2}{\sigma^2} a(\etree 2 1 1 3)+\frac{4}{\sigma^4}a(\etree 1 1 2 1)\rpar \etree 3 2 0 1.
\end{align*}
For order 3, the expression of $\widetilde{\gamma_2}$ can be found in the appendix.
\end{theorem}

\begin{remark}
\label{remark:natural_conditions}
In contrast to Section \ref{Section:IPP}, we choose in Theorem \ref{Theorem:Formal_IPP} not to reduce the tree $\etree 3 2 0 1$ for $\gamma_1$ and the three trees $\etree 4 4 0 1$, $\etree 3 2 1 1$ and $\etree 4 3 0 1$ for $\gamma_2$ in forests with exactly one edge linked to the root. The reason is these trees do not simplify well, and the coefficient multiplying them vanishes for most methods (for all consistent Runge-Kutta methods for example). If one wants to compute the $f_i$ of Section \ref{Section:High_Order_Integrators}, one should integrate by parts these trees first.
\end{remark}

\section{Construction of high order integrators}
\label{section:construction_high_order_integrators}

\subsection{Improvement of a method order via a modified equation}
\label{Section:High_Order_Integrators}
In \cite{Abdulle14hon}, a recursive method to obtain integrators of any order for the convergence to the invariant measure is presented.
Let us suppose we have an integrator of order exactly $p\geq 1$ for the invariant measure, then, using Theorem \ref{thm:orderAstar}, for all $j<p$, $\AA_j^*\rho_\infty=0$ and $\AA_p^*\rho_\infty\neq 0$.
By integrating by parts, we can write $\int_{\R^d} \AA_p \phi\rho_\infty dx$ as $\int_{\R^d} \phi'f_p\rho_\infty dx$. We then consider the same numerical integrator but for the modified equation where we replaced $f$ by $f-h^pf_p$.
Applying Theorem \ref{thm:orderAstar} to the new context, we see that this integrator is at least of order $p+1$ for the original equation.

In this section, we give tools to simplify the computation of those modified integrators, in particular to calculate simply the operators $\AA_j$, and to find the function $f_p$.

\begin{ex*}
For the $\theta$-method \eqref{equation:theta_method}, we have $\AA_1=F(\gamma)$ where $\gamma$ is given by \eqref{equation:gamma_theta_method}.
Applying integration by parts as described in Section \ref{Section:IPP}, we obtain
\begin{align*}
\gamma &\sim \left(\theta-\frac{1}{2}\right) \etree 3 1 0 1
-\frac{1}{2} \gtree 2 1 0 1 \: \etree 2 1 0 1
-\frac{1}{2} \etree 1 2 0 1 \: \etree 2 1 0 1
+\frac{\sigma^2}{2}(1-\theta)\etree 2 1 1 1
+\sigma^2(1-\theta) \gtree 3 3 1 1
+\frac{\sigma^2}{2} \gtree 1 1 1 1 \: \etree 2 1 0 1\\
&+\frac{\sigma^2}{2} \gtree 2 2 1 1 \: \etree 2 1 0 1
-\frac{\sigma^4}{8} \gtree 2 3 2 1
- \frac{\sigma^4}{8} \gtree 1 1 1 1 \: \gtree 2 4 1 1
- \frac{\sigma^4}{4} \gtree 3 2 2 1
-\frac{\sigma^4}{8} \gtree 2 2 1 1 \: \gtree 2 4 1 1.
\end{align*}
Then $f_1$ is given by
\begin{align*}
(f_1)_i&=\left(\theta-\frac{1}{2}\right) f_i'f
-\frac{1}{2} g'f  f_i
-\frac{1}{2} \Div(f) f_i
+\frac{\sigma^2}{2}(1-\theta) \Delta f_i
+\sigma^2(1-\theta) \sum_{j=1}^d \partial_j f_i   \partial_j g\\
&+\frac{\sigma^2}{2} \Delta g   f_i
+\frac{\sigma^2}{2} \sum_{j=1}^d (\partial_j g)^2 f_i
-\frac{\sigma^4}{8}\partial_i \Delta g
-\frac{\sigma^4}{8}\Delta g \partial_i g
-\frac{\sigma^4}{4}\sum_{j=1}^d \partial_j g \partial_{i,j} g\\
&-\frac{\sigma^4}{8}\sum_{j=1}^d (\partial_j g)^2 \partial_i g
.
\end{align*}
\end{ex*}

\begin{remark}
Note that the above computation does not suppose Assumption \ref{Assumption:f_gradient} to be satisfied. But as the function $g$ is not explicitly known in the general case in practice, we cannot implement easily the corresponding scheme.
Under Assumption \ref{Assumption:f_gradient}, we can directly apply Theorem \ref{Theorem:Formal_IPP} and find a simple expression of $f_1$:
$$f_1=\left(\frac{1}{2}-\theta\right)(f'f+\frac{\sigma^2}{2}\Delta f).$$
\end{remark}


\begin{app*}
Let us calculate the modified integrator $f_2$, where $f$ follows Assumption \ref{Assumption:f_gradient}.
First we rewrite the differential operator $\AA_2$ for the modified equation where we replaced $f$ by $f-hf_1$. We call it $\AA_2^{(1)}$.
We find $\AA_2^{(1)}=F(\gamma)$ where
\begin{align*}
\gamma &=
\theta(3\theta-1)\etree 4 1 0 1
+\theta(3\theta-1)\frac{\sigma^2}{2}\etree 3 1 1 1
+\theta(4\theta-1)\frac{\sigma^2}{2}\etree 3 1 1 2
+\theta^2\etree 4 2 0 1
+\theta^2\frac{\sigma^4}{4}\etree 2 1 2 5\\
&+\theta^2\sigma^2\etree 3 1 1 3
+\theta(4\theta-1)\frac{\sigma^2}{2}\etree 3 1 1 4
+\theta(2\theta+1)\frac{\sigma^2}{2}\etree 3 1 1 5
+\theta(2\theta+1)\frac{\sigma^4}{4}\etree 2 1 2 4
+\frac{4\theta-1}{2}\etree 4 3 0 1\\
&+(4\theta-1)\frac{\sigma^2}{4}\etree 3 2 1 1
+\theta^2\frac{\sigma^2}{2}\etree 3 2 1 2
+(4\theta-1)\frac{\sigma^2}{4}\etree 3 1 1 6
+(4\theta-1)\frac{\sigma^4}{8}\etree 2 1 2 1
+\theta\frac{\sigma^4}{2}\etree 2 1 2 2\\
&+\theta\sigma^2\etree 3 2 1 3
+\frac{1}{6}\etree 4 4 0 1
+\theta\frac{\sigma^4}{2}\etree 2 1 2 6
+\frac{\sigma^2}{4}\etree 3 2 1 4
+\frac{\sigma^4}{8}\etree 2 1 2 7
+\frac{\sigma^6}{48}\etree 1 1 3 1
.
\end{align*}
Using Theorem \ref{Theorem:Formal_IPP}, we find 
\begin{align*}
\gamma &\sim \left(-2\theta^2+2\theta-\frac{1}{2}\right)\etree 4 1 0 1
+\left(-\theta^2+\theta-\frac{1}{4}\right)\sigma^2\etree 3 1 1 1
+\left(-\frac{3\theta^2}{2}+\frac{3\theta}{2}-\frac{1}{3}\right)\sigma^2\etree 3 1 1 2\\
&+\left(-\theta^2+\theta-\frac{1}{6}\right)\etree 4 2 0 1
+\left(-\frac{\theta^2}{4}+\frac{\theta}{4}-\frac{1}{24}\right)\sigma^4\etree 2 1 2 5
+\left(-\theta^2+\theta-\frac{1}{6}\right)\sigma^2\etree 3 1 1 3.
\end{align*}
Thus we define
\begin{align*}
f_2&=\left(-2\theta^2+2\theta-\frac{1}{2}\right) f'f'f
+\left(-\theta^2+\theta-\frac{1}{4}\right)\sigma^2 f'\Delta f\\
&+\left(-\frac{3\theta^2}{2}+\frac{3\theta}{2}-\frac{1}{3}\right)\sigma^2 \sum \limits_i f''(e_i,f'(e_i))
+\left(-\theta^2+\theta-\frac{1}{6}\right) f''(f,f)\\
&+\left(-\frac{\theta^2}{4}+\frac{\theta}{4}-\frac{1}{24}\right)\sigma^4 \Delta^2 f
+\left(-\theta^2+\theta-\frac{1}{6}\right)\sigma^2 (\Delta f)'(f),
\end{align*}
and, if the $\theta$-scheme applied to $dX=(f-hf_1-h^2f_2)dt+\sigma dW$ is ergodic, then it has order 3 for the invariant measure.
\end{app*}

For $f_1$, we recover the formula of \cite[Prop.\ts 5.1]{Abdulle14hon}.
The computation of $f_2$ was first done for $\theta=0$ in \cite[Prop.\ts 5.2]{Abdulle14hon}, which reveals a typographical error.

This method can give numerical integrators of any order, but it comes with a high computing price if the partial derivatives of $f$ are difficult to compute. In the following sections, we present order conditions for certain classes of numerical schemes, in order to obtain high order methods avoiding derivatives and unnecessary evaluations of $f$.

\subsection{Order conditions for stochastic Runge-Kutta schemes}
\label{section:Runge_Kutta}
We consider stochastic Runge-Kutta schemes \eqref{equation:defRK} for the overdamped Langevin equation \eqref{equation:sde}.
We set $c_i=\sum \limits_{j=1}^s a_{ij}$.
In this section, we also assume Assumption \ref{Assumption:f_gradient} to simplify the computation.
Using the proposed framework, our goal is to find algebraic conditions on the coefficients $A=(a_{ij})$, $b=(b_i)$ and $d=(d_i)$ to achieve a given order condition for the invariant measure.

First, we suppose $\sum b_i=1$ in order for $\AA_0=\LL$ in Assumption \ref{Assumption:A0=L} to be satisfied.
Then $\AA_1\phi=F(\gamma_1)(\phi)$ where
$$\gamma_1= \sum b_i c_i \etree 3 1 0 1
+\frac{\sigma^2}{2} \sum b_i d_i^2 \etree 2 1 1 1
+\sigma^2 \sum b_i d_i \etree 2 1 1 2
+\frac{1}{2} \etree 3 2 0 1
+\frac{\sigma^2}{2} \etree 2 1 1 3
+\frac{\sigma^4}{8} \etree 1 1 2 1.$$
Theorem \ref{Theorem:Formal_IPP} yields
$$\gamma_1 \sim \left(\sum b_i c_i+\frac{1}{2}-2\sum b_i d_i\right) \etree 3 1 0 1
+\frac{\sigma^2}{2} \left(\sum b_i d_i^2+\frac{1}{2}-2\sum b_i d_i\right)  \etree 2 1 1 1.$$
Thus if we suppose
$$
\sum b_i=1,\quad
\sum b_i c_i+\frac{1}{2}-2\sum b_i d_i=0,\quad
\sum b_i d_i^2+\frac{1}{2}-2\sum b_i d_i=0,
$$
then $\AA_1\phi \sim 0$. We have a Runge-Kutta scheme of order 2.

By continuing this methodology, we obtain the order conditions of order 3, and our analysis allows us to obtain the conditions for any order.
The following theorem states the order conditions for Runge-Kutta methods.

\begin{theorem}
\label{Theorem:RK_order_conditions}
Assume Assumption \ref{Assumption:f_gradient} and consider an ergodic Runge-Kutta method \eqref{equation:defRK} with $\sum b_i=1$. Using the same notation as in Theorem \ref{Theorem:Formal_IPP}, if $A$, $b$ and $d$ are chosen such that $F(\widetilde{\gamma_i})=0$ for all $1\leq i<p$, then the method has at least order $p$ for the invariant measure.
In particular, Table \ref{table:RK_order_conditions} gives sufficient conditions to have consistency and order 2 or 3 for the invariant measure for Runge-Kutta schemes.

\begin{table}[tbh]
\setcellgapes{3pt}
\begin{center}
\begin{tabular}{|c|c|c|l|}
  \hline
  Order & Tree $\tau$ & $F(\tau)(\phi)$ & \multicolumn{1}{c|}{Order condition} \\
  \hline
  1 & $\etree 2 1 0 1$ & $\phi'f$ & $\sum b_i=1$ \\
  \hline
  2 & $\etree 3 1 0 1$ & $\phi'f'f$ & $\sum b_i c_i-2\sum b_i d_i=-\frac{1}{2}$ \\
  \cline{2-4}
   & $\etree 2 1 1 1$ & $\phi'\Delta f$ & $\sum b_i d_i^2-2\sum b_i d_i=-\frac{1}{2}$ \\
  \hline
  3 & $\etree 4 1 0 1$ & $\phi'f'f'f$ & $\sum b_i a_{ij}c_j-2\sum b_i a_{ij}d_j+\sum b_i c_i-\left(\sum b_i d_i\right)^2=0$ \\
  \cline{2-4}
   & $\etree 3 1 1 1$ & $\phi'f'\Delta f$ & $\sum b_i a_{ij}d_j^2-2\sum b_i a_{ij}d_j+\sum b_i c_i-\left(\sum b_i d_i\right)^2=0$ \\
   \cline{2-4}
   &  &  & $\frac{1}{2}\sum b_i c_i^2 -2\sum b_i d_ic_i$ \\[-2ex]
   & $\etree 4 2 0 1$ & $\phi'f''(f,f)$ & $-2\sum b_i d_i+2\sum b_i d_i^2 + \sum b_i c_i=-\frac{1}{3}$ \\
   \cline{2-4}
   &  &  & $\sum b_i d_i a_{ij}d_j-\sum b_i c_id_i-\sum b_i d_i+\sum b_i d_i^2$ \\[-2ex]
   & $\etree 3 1 1 2$ & $\sum \phi'f''(f'(e_i),e_i)$ & $+\sum b_i c_i-\sum b_i a_{ij}d_j-\frac{1}{2}\left(\sum b_i d_i\right)^2=-\frac{1}{6}$ \\
   \cline{2-4}
   &  &  & $\frac{1}{2}\sum b_i c_id_i^2-\sum b_id_i^3-2\sum b_i d_i$ \\[-2ex]
   & $\etree 3 1 1 3$ & $\phi'((\Delta f)'(f))$ & $+\frac{5}{2}\sum b_i d_i^2-\sum b_i c_id_i+\frac{1}{2}\sum b_i c_i=-\frac{1}{3}$ \\
   \cline{2-4}
   &  &  & $\frac{1}{8}\sum b_i d_i^4-\frac{1}{2}\sum b_id_i^3$ \\[-0.5ex]
   & $\etree 2 1 2 5$ & $\phi'\Delta^2 f$ & $-\frac{1}{2}\sum b_i d_i+\frac{3}{4}\sum b_i d_i^2=-\frac{1}{12}$ \\
  \hline
\end{tabular}
\caption{Runge-Kutta order conditions for the invariant measure (See Theorem \ref{Theorem:RK_order_conditions}). The sums are over all involved indices.}
\label{table:RK_order_conditions}
\end{center}
\setcellgapes{1pt}
\end{table}
\end{theorem}

\begin{remark}
To check if a scheme has weak order $p$, one can develop $\LL^i\phi$, $i\leq p$, in exotic aromatic B-series (a simple method for this computation is proposed in Section \ref{postprocessed integrators}) and prove that the forest coefficients of $\frac{1}{i!}\LL^i\phi$ and $\AA_{i-1}\phi$ are equal, yielding the order $p$ estimate \eqref{equation:weak_condition}.
In Table \ref{table:RK_weak_order_conditions}, we collect the corresponding order conditions up to order $p\leq 3$. We recover exactly the same order conditions as first derived in \cite[Thm.\ts 4]{Debrabant10rkm} using different types of trees and B-series.
We recall that the conditions for weak order 3 have no solution for a method of the form \eqref{equation:defRK} with only $l=1$ noise. Indeed, fixing $d^{(2)}=0$ in Table \ref{table:RK_weak_order_conditions}, we obtain the incompatible order conditions $\sum b_i \big(d_i^{(1)}\big)^2=\frac{1}{2}$ (third line of Table \ref{table:RK_weak_order_conditions}) and $\sum b_i \big(d_i^{(1)}\big)^2=\frac{1}{3}$ (last line of Table \ref{table:RK_weak_order_conditions}).
Taking $l=2$ noises in the method \eqref{equation:defRK} is sufficient for reaching weak order 3 for general $f$ satisfying Assumption \ref{Assumption:f_gradient}.

\begin{table}
\setcellgapes{3pt}
\begin{center}
\begin{tabular}{|c|c|c|l|}
  \hline
  Order & Tree $\tau$ & $F(\tau)(\phi)$ & \multicolumn{1}{c|}{Order condition} \\
  \hline
  1 & $\etree 2 1 0 1$ & $\phi'f$ & $\sum b_i=1$ \\
  \hline
  2 & $\etree 3 1 0 1$ & $\phi'f'f$ & $\sum b_i c_i=\frac{1}{2}$ \\
  \cline{2-4}
   & $\etree 2 1 1 1$ & $\phi'\Delta f$ & $\sum b_i \big( d_i^{(1)}\big)^2+\sum b_i \big( d_i^{(2)}\big)^2=\frac{1}{2}$ \\
  \cline{2-4}
   & $\etree 2 1 1 2$ & $\sum\phi''(e_i,f'(e_i))$ & $\sum b_i d_i^{(1)}=\frac{1}{2}$ \\
  \hline
   3 & $\etree 4 1 0 1$ & $\phi'f'f'f$ & $\sum b_i a_{ij}c_j=\frac{1}{6}$ \\
  \cline{2-4}
   & $\etree 3 1 1 1$ & $\phi'f'\Delta f$ & $\sum b_i a_{ij}\big( d_j^{(1)}\big)^2+\sum b_i a_{ij}\big( d_j^{(2)}\big)^2=\frac{1}{6}$ \\
   \cline{2-4}
   & $\etree 4 2 0 1$ & $\phi'f''(f,f)$ & $\sum b_i c_i^2=\frac{1}{3}$ \\
   \cline{2-4}
   & $\etree 3 1 1 2$ & $\sum \phi'f''(f'(e_i),e_i)$ & $\sum b_i d_i^{(1)} a_{ij}d_j^{(1)}+\sum b_i d_i^{(2)} a_{ij}d_j^{(2)}=\frac{1}{6}$ \\
   \cline{2-4}
   & $\etree 3 1 1 3$ & $\phi'((\Delta f)'(f))$ & $\sum b_i c_i \big( d_i^{(1)}\big)^2+\sum b_i c_i\big( d_i^{(2)}\big)^2=\frac{1}{3}$ \\
   \cline{2-4}
   & $\etree 2 1 2 5$ & $\phi'\Delta^2 f$ & $\sum b_i \Big( \big( d_i^{(1)}\big)^2+\big( d_i^{(2)}\big)^2\Big)^2=\frac{1}{3}$ \\
   \cline{2-4}
   & $\etree 3 1 1 4$ & $\sum \phi''(e_i,f'f'(e_i))$ & $\sum b_i a_{ij}d_j^{(1)}=\frac{1}{6}$ \\
   \cline{2-4}
   & $\etree 3 2 1 2$ & $\sum \phi''(f'(e_i),f'(e_i))$ & $\left( \sum b_i d_i^{(2)}\right)^2=\frac{1}{12}$ \\
   \cline{2-4}
   & $\etree 3 1 1 5$ & $\sum \phi''(e_i,f''(e_i,f))$ & $\sum b_i c_i d_i^{(1)}=\frac{1}{3}$ \\
   \cline{2-4}
   & $\etree 2 1 2 4$ & $\sum \phi''(e_i,(\Delta f)'(e_i))$ & $\sum b_i \big( d_i^{(1)}\big)^3+\sum b_i d_i^{(1)} \big( d_i^{(2)}\big)^2=\frac{1}{3}$ \\
   \cline{2-4}
   & $\etree 2 1 2 2$ & $\sum \phi'''(e_i,e_j,f''(e_i,e_j))$ & $\sum b_i \big( d_i^{(1)}\big)^2=\frac{1}{3}$ \\
  \hline
\end{tabular}
\caption{Runge-Kutta standard weak order conditions for $l=2$ noises. The sums are over all involved indices. We recover the same conditions as in \cite{Debrabant10rkm}.}.
\label{table:RK_weak_order_conditions}
\end{center}
\setcellgapes{1pt}
\end{table}
\end{remark}

\begin{remark}
Notice that Assumption \ref{Assumption:f_gradient} permits us to identify the differentials $F(\etree 3 1 1 4)=F(\etree 3 2 1 2)$ and the corresponding order conditions, using Proposition \ref{Proposition:Simplification_rule}. Thus, under Assumption \ref{Assumption:f_gradient}, we can replace the two conditions $\sum b_i a_{ij}d_j=\frac{1}{6}$ and $\left(\sum b_i \widetilde{d_i}\right)^2=\frac{1}{12}$ by one new condition $\sum b_i a_{ij}d_j^{(1)}+\frac{1}{2}\left( \sum b_i d_i^{(2)}\right)^2=\frac{5}{24}$.
\end{remark}

\begin{remark}
As explained in the introduction, the study of weak order conditions using rooted trees is already well documented in the literature, but the framework of exotic aromatic B-series has the advantage to involve rooted forests that do not depend on the dimension~$d$, which permits us to compute integration by parts and hence derive the order conditions for the invariant measure.
Since weak convergence implies convergence with at least the same order for the invariant measure, the weak order conditions (Table \ref{table:RK_weak_order_conditions}) imply the order conditions for the invariant measure (Table \ref{table:RK_order_conditions}).
In particular, comparing Table \ref{table:RK_order_conditions} and Table \ref{table:RK_weak_order_conditions}, we observe that 
there is a lower number of order conditions for the convergence to the invariant measure compared to the standard weak convergence.
\end{remark}

\subsection{Order conditions for postprocessed integrators}
\label{postprocessed integrators}
In this section, we extend our analysis to the case of integrators combined with postprocessors \cite{Vilmart15pif}. 
As stated in Theorem \ref{Theorem:Postprocessed_scheme_order_condition},
it permits us to increase the order for the invariant measure of a given method while maintaining a low number of function evaluations per time step.
We show that exotic aromatic B-series simplify this approach, but one issue remains: the computation of the Lie bracket $[\LL,\overline{\AA_p}]\phi$. 
This is done by the following theorem for the composition of exotic aromatic forests and based on the Leibniz rule.

\begin{theorem}
\label{Theorem:Forest_composition}
Let $\gamma_1$ and $\gamma_2$ be two exotic aromatic rooted forests, and let $\phi: \R^d \rightarrow \R$ be a smooth function. For $\varphi: \pi(r_2)\to V_1$ and $\psi: \Gamma(r_2)\to V_1$, we build $\gamma_{\varphi,\psi}$ by plugging all the edges connected to $r_2$ to the nodes of $\gamma_1$ according to $\varphi$, and all the lianas (counting multiplicity) according to $\psi$. Then the composition of forests is given by
$$F(\gamma_2)(F(\gamma_1)(\phi))=\sum_{\underset{\psi: \Gamma(r_2)\to V_1}{\varphi: \pi(r_2)\to V_1}} F(\gamma_{\varphi,\psi})(\phi).$$
\end{theorem}

Various composition rules for B-series and aromatic B-series have been studied in the literature (see \cite{Hairer06gni,Chartier10aso,Bogfjellmo15aso} and the references therein). The main difference from these previous works is that we compose only the roots of exotic aromatic rooted forests, because this corresponds to composing linear differential operators.

\begin{proof}
Using Definition \ref{elementary_differential}, we have
$$F(\gamma_1)(\phi)=\sum_{i_{v_1^{(1)}},\dots,i_{v_{m_1}^{(1)}}} \sum_{j_{l_1^{(1)}},\dots,j_{l_{s_1}^{(1)}}}
\left(\prod_{v\in V_1^0} \partial_{I_{\pi(v)}} \partial_{J_{\Gamma(v)}} f_{i_v}\right)
\partial_{I_{\pi(r_1)}} \partial_{J_{\Gamma(r_1)}} \phi.$$
Then we replace $\phi$ by $F(\gamma_1)(\phi)$ and use the Leibniz rule to distribute the partial derivatives.
\begin{align*}
F(\gamma_2)(F(\gamma_1)(\phi))&=
\sum_{\underset{i_{v_1^{(2)}},\dots,i_{v_{m_2}^{(2)}}}{i_{v_1^{(1)}},\dots,i_{v_{m_1}^{(1)}}}}
\sum_{\underset{j_{l_1^{(2)}},\dots,j_{l_{s_2}^{(2)}}}{j_{l_1^{(1)}},\dots,j_{l_{s_1}^{(1)}}}}
\left(\prod_{v\in V_2^0} \partial_{I_{\pi(v)}} \partial_{J_{\Gamma(v)}} f_{i_v}\right)\\
&\cdot\partial_{I_{\pi(r_2)}} \partial_{J_{\Gamma(r_2)}}\left[\left(\prod_{v\in V_1^0} \partial_{I_{\pi(v)}} \partial_{J_{\Gamma(v)}} f_{i_v}\right) \partial_{I_{\pi(r_1)}} \partial_{J_{\Gamma(r_1)}} \phi\right]\\
&=\sum_{\underset{\psi: \Gamma(r_2)\to V_1}{\varphi: \pi(r_2)\to V_1}}
\sum_{\underset{i_{v_1^{(2)}},\dots,i_{v_{m_2}^{(2)}}}{i_{v_1^{(1)}},\dots,i_{v_{m_1}^{(1)}}}}
\sum_{\underset{j_{l_1^{(2)}},\dots,j_{l_{s_2}^{(2)}}}{j_{l_1^{(1)}},\dots,j_{l_{s_1}^{(1)}}}}
\left(\prod_{v\in V_2^0} \partial_{I_{\pi(v)}} \partial_{J_{\Gamma(v)}} f_{i_v}\right)\\
&\cdot \left(\prod_{v\in V_1^0} \partial_{I_{\pi(v)\cup \varphi^{-1}(\{v\})}} \partial_{J_{\Gamma(v)\cup \psi^{-1}(\{v\})}} f_{i_v}\right) \partial_{I_{\pi(r_1)\cup \varphi^{-1}(\{r_1\})}} \partial_{J_{\Gamma(r_1)\cup \psi^{-1}(\{r_1\})}} \phi\\
&=\sum_{\underset{\psi: \Gamma(r_2)\to V_1}{\varphi: \pi(r_2)\to V_1}} F(\gamma_{\varphi,\psi})(\phi).
\end{align*}
\vskip-6ex
\end{proof}

\begin{ex*}
We recall $\LL=F(\etree 2 1 0 1+\frac{\sigma^2}{2}\etree 1 1 1 1)$, then we can compute $\LL^2\phi$.
Using Theorem \ref{Theorem:Forest_composition}, we obtain
$$
F(\etree 2 1 0 1)(\LL \phi)=F(\etree 3 1 0 1+\etree 3 2 0 1 +\frac{\sigma^2}{2}\etree 2 1 1 3)(\phi)
$$
and
$$
F(\etree 1 1 1 1)(\LL \phi)=F(\etree 2 1 1 3+\etree 2 1 1 1 +2\etree 2 1 1 2+\frac{\sigma^2}{2}\etree 1 1 2 1)(\phi).
$$
Combining the two previous equalities, we deduce
$$
\LL^2\phi=F(\etree 3 1 0 1+\etree 3 2 0 1 +\sigma^2\etree 2 1 1 3+\frac{\sigma^2}{2}\etree 2 1 1 1 +\sigma^2\etree 2 1 1 2+\frac{\sigma^4}{4}\etree 1 1 2 1)(\phi).
$$
\end{ex*}

Using Theorems \ref{Theorem:Postprocessed_scheme_order_condition} and \ref{Theorem:Forest_composition}, we obtain general conditions on postprocessors to increase by 1 the order of a given method.

\begin{theorem}
\label{Theorem:General_order_conditions_with_postprocessor}
Using notation and assumptions of Theorem \ref{Theorem:Postprocessed_scheme_order_condition}, if the numerical scheme and the postprocessor can be developed in exotic aromatic B-series of the respective forms \ref{equation:dvp_U_Bseries} and
$$
\E[\phi(G_n(x))]=F(\etree 1 1 0 1)(\phi) + \sum_{\underset{1\leq\abs{\gamma}\leq p}{\gamma \in \EE\AA\TT}} h^{\abs{\gamma}} \overline{a}(\gamma) F(\gamma)(\phi) + \dots,
$$
if we denote $\gamma$ the exotic aromatic B-series such that $F(\gamma)=(\AA_p+[\LL,\overline{\AA_p}])$ and if $\gamma\sim 0$,
then $\overline{X_n}$ is of order $p+1$ for the invariant measure.
In particular, if the order 2 conditions in Table \ref{table:General_order_2_condition_with_postprocessor} are verified, then the method has order 2 for the invariant measure.
If we suppose the order 2 for the invariant measure of the numerical method and $\overline{a}(\etree 2 1 0 1)+\frac{2}{\sigma^2}\overline{a}(\etree 1 1 1 1)=0$ (in order to have $\overline{\AA_0}=\alpha \LL$), and if conditions in Table \ref{table:General_order_3_condition_with_postprocessor} are verified, then the method is of order 3 for the invariant measure.

\begin{table}[tbh]
\setcellgapes{3pt}
\begin{center}
\begin{tabular}{|c|l|}
  \hline
Tree $\gamma$ & \multicolumn{1}{c|}{Order condition} \\
  \hline
  $\etree 3 1 0 1$ & $a(\etree 3 1 0 1)-\frac{2}{\sigma^2}a(\etree 2 1 1 2)+\frac{2}{\sigma^2}a(\etree 2 1 1 3)-\frac{4}{\sigma^4}a(\etree 1 1 2 1)-2\overline{a}(\etree 2 1 0 1)+\frac{4}{\sigma^2}\overline{a}(\etree 1 1 1 1)=0$ \\
  \hline
  $\etree 2 1 1 1$ & $a(\etree 2 1 1 1)-a(\etree 2 1 1 2)+a(\etree 2 1 1 3)-\frac{2}{\sigma^2}a(\etree 1 1 2 1)-\frac{\sigma^2}{2}\overline{a}(\etree 2 1 0 1)+\overline{a}(\etree 1 1 1 1)=0$ \\
  \hline
  $\etree 3 2 0 1$ & $a(\etree 3 2 0 1)-\frac{2}{\sigma^2} a(\etree 2 1 1 3)+\frac{4}{\sigma^4}a(\etree 1 1 2 1)=0$ \\
  \hline
\end{tabular}
\caption{General order 2 conditions with postprocessor (See Theorem \ref{Theorem:General_order_conditions_with_postprocessor}).}
\label{table:General_order_2_condition_with_postprocessor}
\end{center}
\setcellgapes{1pt}
\end{table}

\begin{table}
\setcellgapes{3pt}
\begin{center}
\begin{tabular}{|c|l|}
  \hline
  Tree $\tau$ & \multicolumn{1}{c|}{Order condition} \\
  \hline
   & $a(\etree 4 1 0 1)-\frac{2}{\sigma^2}a(\etree 3 1 1 4)-\frac{2}{\sigma^2}a(\etree 3 2 1 2)+\frac{2}{\sigma^2}a(\etree 3 1 1 6)+\frac{2}{\sigma^2}a(\etree 3 2 1 3)-\frac{4}{\sigma^2}a(\etree 3 2 1 4)-\frac{4}{\sigma^4}a(\etree 2 1 2 6)$ \\[-1.5ex]
  $\etree 4 1 0 1$ & $+\frac{12}{\sigma^4}a(\etree 2 1 2 7)-\frac{24}{\sigma^6}a(\etree 1 1 3 1)-2\overline{a}(\etree 3 1 0 1)+\frac{4}{\sigma^2}\overline{a}(\etree 2 1 1 2)+2\overline{a}(\etree 3 2 0 1)-\frac{8}{\sigma^2}\overline{a}(\etree 2 1 1 3)+\frac{16}{\sigma^4}\overline{a}(\etree 1 1 2 1)=0$ \\
  \hline
   & $a(\etree 3 1 1 1)-a(\etree 3 1 1 4)-a(\etree 3 2 1 2)+a(\etree 3 1 1 6)+a(\etree 3 2 1 3)-2 a(\etree 3 2 1 4)-\frac{2}{\sigma^2}a(\etree 2 1 2 6)+\frac{6}{\sigma^2}a(\etree 2 1 2 7)$ \\[-1ex]
  $\etree 3 1 1 1$ & $-\frac{12}{\sigma^4}a(\etree 1 1 3 1)-\frac{\sigma^2}{2}\overline{a}(\etree 3 1 0 1)-\overline{a}(\etree 2 1 1 1)+2\overline{a}(\etree 2 1 1 2)+\sigma^2\overline{a}(\etree 3 2 0 1)-4\overline{a}(\etree 2 1 1 3)+\frac{8}{\sigma^2}\overline{a}(\etree 1 1 2 1)=0$ \\
  \hline
   & $a(\etree 4 2 0 1)-\frac{2}{\sigma^2} a(\etree 3 1 1 5)+\frac{4}{\sigma^4} a(\etree 2 1 2 2)+\frac{2}{\sigma^2}a(\etree 3 1 1 6)-\frac{4}{\sigma^4} a(\etree 2 1 2 6)$ \\
  $\etree 4 2 0 1$ & $+\frac{4}{\sigma^4} a(\etree 2 1 2 7)-\frac{8}{\sigma^6}a(\etree 1 1 3 1)-\overline{a}(\etree 3 1 0 1)-\overline{a}(\etree 3 2 0 1)+\frac{2}{\sigma^2}\overline{a}(\etree 2 1 1 2)=0$ \\
  \hline
   & $a(\etree 3 1 1 2)-a(\etree 3 1 1 5)-a(\etree 3 2 1 2)-a(\etree 3 1 1 4)+2a(\etree 3 1 1 6)+\frac{2}{\sigma^2} a(\etree 2 1 2 2)+a(\etree 3 2 1 3)-2a(\etree 3 2 1 4)-\frac{4}{\sigma^2}a(\etree 2 1 2 6)$ \\
  $\etree 3 1 1 2$ & $+\frac{8}{\sigma^2} a(\etree 2 1 2 7)-\frac{16}{\sigma^4} a(\etree 1 1 3 1)-\sigma^2\overline{a}(\etree 3 1 0 1)+\sigma^2\overline{a}(\etree 3 2 0 1)+2\overline{a}(\etree 2 1 1 2)-4\overline{a}(\etree 2 1 1 3)+\frac{8}{\sigma^2}\overline{a}(\etree 1 1 2 1)=0$ \\
  \hline
   & $a(\etree 3 1 1 3)-\frac{2}{\sigma^2}a(\etree 2 1 2 4)-a(\etree 3 1 1 5)+\frac{2}{\sigma^2}a(\etree 2 1 2 1)+a(\etree 3 1 1 6)+\frac{4}{\sigma^2} a(\etree 2 1 2 2)-\frac{4}{\sigma^2}a(\etree 2 1 2 6)$ \\
  $\etree 3 1 1 3$ & $+\frac{4}{\sigma^2} a(\etree 2 1 2 7)-\frac{8}{\sigma^4} a(\etree 1 1 3 1)-\frac{\sigma^2}{2}\overline{a}(\etree 3 1 0 1)-\overline{a}(\etree 2 1 1 1)+2\overline{a}(\etree 2 1 1 2)-2\overline{a}(\etree 2 1 1 3)+\frac{4}{\sigma^2}\overline{a}(\etree 1 1 2 1)=0$ \\
  \hline
   & $a(\etree 2 1 2 5)-a(\etree 2 1 2 4)+a(\etree 2 1 2 1)+ a(\etree 2 1 2 2)-a(\etree 2 1 2 6)+ a(\etree 2 1 2 7)$ \\
  $\etree 2 1 2 5$ & $-\frac{2}{\sigma^2} a(\etree 1 1 3 1)-\frac{\sigma^2}{2}\overline{a}(\etree 2 1 1 1)+\frac{\sigma^2}{2}\overline{a}(\etree 2 1 1 2)-\frac{\sigma^2}{2}\overline{a}(\etree 2 1 1 3)+\overline{a}(\etree 1 1 2 1)=0$ \\
  \hline
  $\etree 4 4 0 1$ & $a(\etree 4 4 0 1)-\frac{2}{\sigma^2}a(\etree 3 2 1 4)+\frac{4}{\sigma^4} a(\etree 2 1 2 7)-\frac{8}{\sigma^6}a(\etree 1 1 3 1)=0$ \\
  \hline
   & $a(\etree 3 2 1 1)-a(\etree 3 2 1 3)+2 a(\etree 3 2 1 4)-\frac{2}{\sigma^2}a(\etree 2 1 2 1)+\frac{2}{\sigma^2}a(\etree 2 1 2 6)$ \\
  $\etree 3 2 1 1$ & $-\frac{6}{\sigma^2}a(\etree 2 1 2 7)+\frac{12}{\sigma^4}a(\etree 1 1 3 1)-\sigma^2\overline{a}(\etree 3 2 0 1)+2\overline{a}(\etree 2 1 1 3)-\frac{4}{\sigma^2}\overline{a}(\etree 1 1 2 1)=0$ \\
  \hline
   & $a(\etree 4 3 0 1)-\frac{2}{\sigma^2}a(\etree 3 2 1 3)-\frac{2}{\sigma^2}a(\etree 3 1 1 6)+\frac{4}{\sigma^4}a(\etree 2 1 2 6)+\frac{4}{\sigma^2} a(\etree 3 2 1 4)$ \\
  $\etree 4 3 0 1$ & $-\frac{12}{\sigma^4}a(\etree 2 1 2 7)+\frac{24}{\sigma^6}a(\etree 1 1 3 1)-4\overline{a}(\etree 3 2 0 1)+\frac{8}{\sigma^2}\overline{a}(\etree 2 1 1 3)-\frac{16}{\sigma^4}\overline{a}(\etree 1 1 2 1)=0$ \\
  \hline
\end{tabular}
\caption{General order 3 conditions with postprocessor (See Theorem \ref{Theorem:General_order_conditions_with_postprocessor}).}
\label{table:General_order_3_condition_with_postprocessor}
\end{center}
\setcellgapes{1pt}
\end{table}
\end{theorem}

\begin{theorem}
\label{Theorem:Postprocessed_Runge_Kutta_order_conditions}
Consider an ergodic Runge-Kutta method of order $p\geq 1$ for the invariant measure of the form \eqref{equation:defRK} and the following associated postprocessor
$$
\begin{array}{l}
  \overline{Y_i}=X_n+h\sum \limits_{j=1}^s \overline{a_{ij}}f(\overline{Y_j}) +\overline{d_i}\sigma\sqrt{h}\,\overline{\xi_n} , \qquad i = 1, \dots,s, \\
  \overline{X_n}=X_n+h\sum \limits_{i=1}^s \overline{b_i} f(\overline{Y_i}) +\overline{d_0}\sigma\sqrt{h}\,\overline{\xi_n}.
\end{array}
$$
Assume Assumption \ref{Assumption:f_gradient}. If $\gamma$ is the exotic aromatic B-series such that $F(\gamma)=(\AA_p+[\LL,\overline{\AA_p}])$ and if $A$, $b$, $d$, $\overline{A}$, $\overline{b}$, $\overline{d}$, $\overline{d_0}$ are chosen such that $\gamma \sim 0$ then the postprocessed method $\overline{X_n}$ has at least order $p+1$ for the invariant measure.
In particular, if the conditions of order 2 in Table \ref{table:Order_conditions_for_Runge-Kutta_method_with_Runge-Kutta_postprocessor} are verified, then the postprocessed integrator has order 2.
If the Runge-Kutta method has order 2 for the invariant measure (see Table \ref{table:RK_order_conditions}), if $\sum \overline{b_i}=\overline{d_0}^2$ and if the conditions of order 3 in Table \ref{table:Order_conditions_for_Runge-Kutta_method_with_Runge-Kutta_postprocessor} are verified, then the method has order 3.

\begin{table}[tbh]
\setcellgapes{3pt}
\begin{center}
\begin{tabular}{|c|c|l|}
  \hline
  Order & Tree $\tau$ & \multicolumn{1}{c|}{Order condition} \\
  \hline
  2 & $\etree 3 1 0 1$ & $\sum b_i c_i-2\sum b_i d_i-2\sum \overline{b_i}+2\overline{d_0}^2=-\frac{1}{2}$ \\
  \cline{2-3}
   & $\etree 2 1 1 1$ & $\sum b_i d_i^2-2\sum b_i d_i-\sum \overline{b_i}+\overline{d_0}^2=-\frac{1}{2}$ \\
  \hline
   &  & $\sum b_i a_{ij}c_j-2\sum b_i a_{ij}d_j+\sum b_i c_i$\\[-3ex]
  3 & $\etree 4 1 0 1$ & $-\left(\sum b_i d_i\right)^2-2\sum \overline{b_i}\overline{c_i}+4\overline{d_0}\sum \overline{b_i}\overline{d_i}-\overline{d_0}^4=0$ \\
  \cline{2-3}
   &  & $\sum b_i a_{ij}d_j^2-2\sum b_i a_{ij}d_j+\sum b_i c_i$\\[-3ex]
   & $\etree 3 1 1 1$ & $-\left(\sum b_i d_i\right)^2-\sum \overline{b_i}\overline{c_i}-\sum \overline{b_i}\overline{d_i}^2+4\overline{d_0}\sum \overline{b_i}\overline{d_i}-\overline{d_0}^4=0$ \\
  \cline{2-3}
   &  & $\frac{1}{2}\sum b_i c_i^2 -2\sum b_i d_ic_i-2\sum b_i d_i+2\sum b_i d_i^2 $\\[-2ex]
   & $\etree 4 2 0 1$ & $ + \sum b_i c_i -\sum \overline{b_i}\overline{c_i}+2\overline{d_0}\sum \overline{b_i}\overline{d_i}- \frac{\overline{d_0}^4}{2}=-\frac{1}{3}$ \\
  \cline{2-3}
   &  & $\sum b_i d_i a_{ij}d_j-\sum b_i c_id_i-\sum b_i d_i+\sum b_i d_i^2+\sum b_i c_i$ \\[-2ex]
   & $\etree 3 1 1 2$ & $-\sum b_i a_{ij}d_j-\frac{1}{2}\left(\sum b_i d_i\right)^2-\sum \overline{b_i}\overline{c_i}+2\overline{d_0}\sum \overline{b_i}\overline{d_i}-\frac{\overline{d_0}^4}{2}=-\frac{1}{6}$ \\
  \cline{2-3}
   &  & $\frac{1}{2}\sum b_i c_id_i^2-\sum b_id_i^3-2\sum b_i d_i+\frac{5}{2}\sum b_i d_i^2-\sum b_i c_id_i$ \\[-2ex]
   & $\etree 3 1 1 3$ & $+\frac{1}{2}\sum b_i c_i-\frac{1}{2}\sum \overline{b_i}\overline{c_i}-\frac{1}{2}\sum \overline{b_i}\overline{d_i}^2+2\overline{d_0}\sum \overline{b_i}\overline{d_i}- \frac{\overline{d_0}^4}{2}=-\frac{1}{3}$ \\
  \cline{2-3}
   &  & $\frac{1}{8}\sum b_i d_i^4-\frac{1}{2}\sum b_id_i^3-\frac{1}{2}\sum b_i d_i$\\[-1ex]
   & $\etree 2 1 2 5$ & $+\frac{3}{4}\sum b_i d_i^2-\frac{1}{4}\sum \overline{b_i}\overline{d_i}^2+\frac{1}{2}\overline{d_0}\sum \overline{b_i}\overline{d_i}-\frac{\overline{d_0}^4}{8}=-\frac{1}{12}$ \\
  \hline
\end{tabular}
\caption{Order conditions for Runge-Kutta method with Runge-Kutta postprocessor (See Theorem \ref{Theorem:Postprocessed_Runge_Kutta_order_conditions}). The sums are over all involved indices.}
\label{table:Order_conditions_for_Runge-Kutta_method_with_Runge-Kutta_postprocessor}
\end{center}
\setcellgapes{1pt}
\end{table}
\end{theorem}

\begin{remark}
The condition of $\etree 3 2 0 1$ in Table \ref{table:General_order_2_condition_with_postprocessor} is not modified by the postprocessor as it does not depend on any $\overline{a}(\gamma)$. Thus if the scheme is a consistent Runge-Kutta method, this condition is automatically satisfied (see also Remark \ref{remark:natural_conditions}).
If the postprocessor satisfies the equation
\begin{equation}
\label{equation:formal_pp_condition}
\overline{a}(\etree 3 2 0 1)-\frac{2}{\sigma^2}\overline{a}(\etree 2 1 1 3)+\frac{4}{\sigma^4}\overline{a}(\etree 1 1 2 1)=0,
\end{equation}
then the conditions of $\etree 4 4 0 1$, $\etree 3 2 1 1$ and $\etree 4 3 0 1$ of Table \ref{table:General_order_3_condition_with_postprocessor} are not modified by the postprocessor, and thus are automatically satisfied for consistent Runge-Kutta methods.
Equation \eqref{equation:formal_pp_condition} is verified for the class of Runge-Kutta methods for postprocessors presented in Theorem \ref{Theorem:Postprocessed_Runge_Kutta_order_conditions} that satisfy $\overline{\AA_0}=\alpha \LL$.
\end{remark}

\begin{ex*}
Under Assumption \ref{Assumption:f_gradient}, the following Runge-Kutta method, introduced in \cite{Vilmart15pif}, is of order 2 for the invariant measure of \eqref{equation:Langevin} (if it is ergodic).
$$
\begin{array}{rl}
  X_{n+1}&=X_n+hf(X_{n+1}+\frac{-1+\sqrt{2}}{2}\sigma\sqrt{h}\xi_n)+\sigma\sqrt{h}\xi_n,\\
  \overline{X_n}&=X_n+h\frac{\sqrt{2}}{2}f(\overline{X_n})+\frac{\sqrt{4\sqrt{2}-1}}{2}\sigma\sqrt{h}\,\overline{\xi_n}.
\end{array}
$$
Indeed, its coefficients, placed in the following Butcher tableau, fulfil the conditions of order 2 of Theorem \ref{Theorem:Postprocessed_Runge_Kutta_order_conditions} (See Table \ref{table:Order_conditions_for_Runge-Kutta_method_with_Runge-Kutta_postprocessor}).
$$\begin{array}
{c|c|c|c|c|c}
c & A & d & \overline{c} & \overline{A} & \overline{d}\\
\hline
  & b &   &   & \overline{b} & \overline{d}_0
\end{array}
=
\begin{array}
{c|c|c|c|c|c}
1 & 1 & \frac{1+\sqrt{2}}{2} & \frac{\sqrt{2}}{2} & \frac{\sqrt{2}}{2} & \frac{\sqrt{4\sqrt{2}-1}}{2}\\
\hline
  & 1 &   &   & \frac{\sqrt{2}}{2} & \frac{\sqrt{4\sqrt{2}-1}}{2}
\end{array}$$
\end{ex*}

\subsection{Order conditions for partitioned methods}
\label{section:partitioned_methods}
In \eqref{equation:sde}, we assume $f=f_1+f_2$ and we consider partitioned integrators that apply different numerical treatments to each $f_i$. We explain in this section how to extend the exotic aromatic B-series formalism to compute order conditions for such partitioned integrators.
The advantage is to treat differently each part of $f$ according to their properties. For example, if $f_1$ is stiff and $f_2$ is non-stiff, one would like to apply an implicit method to $f_1$ and an explicit method for $f_2$ (IMEX methods).

Here we follow the formalism of \cite[Sect.\ts III.2]{Hairer06gni} for bicoloured B-series, called P-series. We introduce white nodes ; they represent the function $f_2$. Black nodes now correspond to $f_1$ but the root still corresponds to $\phi$. We call these new forests exotic aromatic P-forests.
There are two slight changes in the computation rules compared to the non-partitioned case:
\begin{itemize}
\item Simplification rule: if $f_1=\nabla V_1$ and $f_2=\nabla V_2$ are both gradients, then
$$\includegraphics[scale=0.5]{Other_trees/simplificationrule5.eps} \sim \frac{2}{\sigma^2} \left(\includegraphics[scale=0.5]{Other_trees/simplificationrule6.eps} + \includegraphics[scale=0.5]{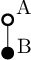}\right).$$
Furthermore, the node $B$ can be replaced by an aromatic root or a white node.
\item The operator $\LL$ is now written as
$$\LL =F(\etree 2 1 0 1+\ptree 2 1 0 1 1+\frac{\sigma^2}{2}\etree 1 1 1 1).$$
\end{itemize}
In addition to the partitioning of the method, one can also add a postprocessor. The results of Section \ref{section:formal_ipp}, \ref{section:Runge_Kutta} and \ref{postprocessed integrators} are straightforwardly adapted to the P-forests.

\begin{theorem}
\label{theorem:RK_partitioned}
Consider a Runge-Kutta method of order $p$ of the form
$$
\begin{array}{rl}
  Y_i&=X_n+h\sum \limits_{j=1}^s a_{ij}f_1(Y_j)+\widehat{a_{ij}}f_2(Y_j) +d_i\sigma\sqrt{h}\xi_n , \qquad i = 1, \dots,s, \\
  X_{n+1}&=X_n+h\sum \limits_{i=1}^s b_i f_1(Y_i)+\widehat{b_i} f_2(Y_i) +\sigma\sqrt{h}\xi_n,
\end{array}
$$
together with the following Runge-Kutta postprocessor
$$
\begin{array}{rl}
  \overline{Y_i}&=X_n+h\sum \limits_{j=1}^s \overline{a_{ij}}f_1(\overline{Y_j})+\widehat{\overline{a_{ij}}}f_2(\overline{Y_j}) +\overline{d_i}\sigma\sqrt{h}\,\overline{\xi_n} , \qquad i = 1, \dots,s, \\
  \overline{X_n}&=X_n+h\sum \limits_{i=1}^s \overline{b_i} f_1(\overline{Y_i})+\widehat{\overline{b_i}} f_2(\overline{Y_i})+\overline{d_0}\sigma\sqrt{h}\,\overline{\xi_n}.
\end{array}
$$
Under the notation and assumptions of Theorem \ref{Theorem:Postprocessed_Runge_Kutta_order_conditions}, if we suppose $f_1$ and $f_2$ are gradients, if we choose $f_1$, $f_2$ and the coefficients of the method such that $\gamma \sim 0$ then the method has at least order $p+1$ for the invariant measure.
In particular, the conditions for consistency and order 2 are in Table \ref{table:Order_conditions_for_partitioned_Runge-Kutta_method_with_postprocessor}.

\begin{table}[tbh]
\setcellgapes{3pt}
\begin{center}
\begin{tabular}{|c|c|c|l|}
  \hline
  Order & Tree $\tau$ & $F(\tau)(\phi)$ & \multicolumn{1}{c|}{Order condition} \\
  \hline
  1 & $\etree 2 1 0 1$ & $\phi'f_1$ & $\sum b_i=1$ \\
  \cline{2-4}
   & $\ptree 2 1 0 1 1$ & $\phi'f_2$ & $\sum \widehat{b_i}=1$ \\
  \hline
  2 & $\etree 3 1 0 1$ & $\phi'f_1'f_1$ & $\sum b_i c_i-2\sum b_i d_i-2\sum \overline{b_i}+2\overline{d_0}^2=-\frac{1}{2}$ \\
  \cline{2-4}
   & $\ptree 3 1 0 1 1$ & $\phi'f_1'f_2$ & $\sum b_i \widehat{c_i}-2\sum b_i d_i-\sum \overline{b_i}-\sum \widehat{\overline{b_i}}+2\overline{d_0}^2=-\frac{1}{2}$ \\
  \cline{2-4}
   & $\ptree 3 1 0 1 2$ & $\phi'f_2'f_1$ & $\sum \widehat{b_i} c_i-2\sum \widehat{b_i} d_i-\sum \overline{b_i}-\sum \widehat{\overline{b_i}}+2\overline{d_0}^2=-\frac{1}{2}$ \\
  \cline{2-4}
   & $\ptree 3 1 0 1 3$ & $\phi'f_2'f_2$ & $\sum \widehat{b_i} \widehat{c_i}-2\sum \widehat{b_i} d_i-2\sum \widehat{\overline{b_i}}+2\overline{d_0}^2=-\frac{1}{2}$ \\
  \cline{2-4}
   & $\etree 2 1 1 1$ & $\phi'\Delta f_1$ & $\sum b_i d_i^2-2\sum b_i d_i-\sum \overline{b_i}+\overline{d_0}^2=-\frac{1}{2}$ \\
  \cline{2-4}
   & $\ptree 2 1 1 1 1$ & $\phi'\Delta f_2$ & $\sum \widehat{b_i} d_i^2-2\sum \widehat{b_i} d_i-\sum \widehat{\overline{b_i}}+\overline{d_0}^2=-\frac{1}{2}$ \\
  \hline
\end{tabular}
\caption{Order conditions for partitioned Runge-Kutta method with postprocessor (See Theorem \ref{theorem:RK_partitioned}). The sums are over all involved indices.}
\label{table:Order_conditions_for_partitioned_Runge-Kutta_method_with_postprocessor}
\end{center}
\setcellgapes{1pt}
\end{table}
\end{theorem}

\begin{ex*}
Using the previously introduced formalism, we see that, if $f_1$ and $f_2$ are gradients and $f=f_1+f_2$ satisfies Assumption \ref{Assumption:f_gradient}, the following method, adapted from \cite[Lemma 2.9]{Brehier16hoi}, is of order 2 for the invariant measure of \eqref{equation:Langevin} (if it is ergodic).
$$
\begin{array}{rl}
  X_{n+1} &=X_n+\frac{h}{2}f_1(X_{n+1}+\frac{1}{2}\sigma\sqrt{h}\xi_n)+\frac{h}{2}f_1(X_{n+1}+\frac{3}{2}\sigma\sqrt{h}\xi_n)\\
  &+h f_2(X_n+\frac{1}{2}\sigma\sqrt{h}\xi_n)+\sigma\sqrt{h}\xi_n,\\
  \overline{X_n} &=X_n+\frac{1}{2}\sigma\sqrt{h}\,\overline{\xi_n}.
\end{array}
$$
It can be put in Runge-Kutta form with the coefficients below:
$$\begin{array}
{c|c|c|c|c}
c & A & \widehat{c} & \widehat{A} & d\\
\hline
  & b &   & \widehat{b} &
\end{array}
=
\begin{array}
{c|ccc|c|ccc|c}
0 & 0 & 0 & 0 & 0 & 0 & 0 & 0 & 1/2 \\
1 & 0 & 1/2 & 1/2 & 1 & 1 & 0 & 0 & 1/2 \\
1 & 0 & 1/2 & 1/2 & 1 & 1 & 0 & 0 & 3/2 \\
\hline
 & 0 & 1/2 & 1/2 &  & 1 & 0 & 0 &
\end{array}$$
and $s=0$ and $\overline{d_0}=\frac{1}{2}$ for the postprocessor.

If we add a family of independent noises $(\chi_n)_n$ independent of $(\xi_n)_n$, then by extending Theorem \ref{theorem:RK_partitioned}, we can show that, under the same hypothesis as the previous example, the following IMEX method has order 2 for the invariant measure of \eqref{equation:Langevin}:
$$
\begin{array}{rl}
  X_{n+1}&=X_n+hf_1(X_{n+1}+\frac{1}{2}\sigma\sqrt{h}\chi_n)+h f_2(X_n+\frac{1}{2}\sigma\sqrt{h}\xi_n)+\sigma\sqrt{h}\xi_n,\\
  \overline{X_n}&=X_n+\frac{1}{2}\sigma\sqrt{h}\,\overline{\xi_n}.
\end{array}
$$
\end{ex*}

\subsection{Non-reversible perturbation}
\label{Section:non_reversible_perturbation}
An interesting modification of \eqref{equation:Langevin} is to introduce a non gradient perturbation that preserves the invariant measure. It permits for some classes of problems to improve the rate of convergence to equilibrium \cite{Lelievre13onr}, and it can also reduce the variance \cite{Duncan16vru}.
As in Section \ref{section:partitioned_methods}, we consider the equation \eqref{equation:sde} where $f=f_1+f_2$ and we use bicoloured forests. We suppose $f_1=-\nabla V$ is a gradient, and $f_2$ is a perturbation of $f_1$ that satisfies
\begin{equation}
\label{equation:f2_perturbed}
\Div \left(f_2 e^{-\frac{2}{\sigma^2}V}\right)=0.
\end{equation}
The perturbation $f_2$ does not modify the invariant measure. Indeed equation \eqref{equation:f2_perturbed} implies that the adjoint of $\BB \phi= \phi'(f_2)$ satisfies $\BB^*\rho_\infty=0$, and thus the invariant measure is preserved. A simple example of such non gradient perturbation is $f_2=J\nabla V$, with $J$ a fixed antisymmetric matrix.
We can now apply all the results of Section \ref{section:using_eat} that do not use Assumption \ref{Assumption:f_gradient}.
We have the following useful properties.
\begin{itemize}
\item We still have the simplification rule (see Proposition \ref{Proposition:Simplification_rule}):
$\includegraphics[scale=0.5]{Other_trees/simplificationrule5.eps} \sim \frac{2}{\sigma^2} \includegraphics[scale=0.5]{Other_trees/simplificationrule6.eps}$.
Furthermore, the node $B$ can be replaced by an aromatic root or a white node.
\item The generator reads
$\LL =F(\etree 2 1 0 1+\ptree 2 1 0 1 1+\frac{\sigma^2}{2}\etree 1 1 1 1)$.
\item We have
$F(\ptree 1 2 0 1 1)=-F(\ptree 2 1 0 1 2)$,
and these differentials vanish if $f_2=J\nabla V$.
\end{itemize}
The two first properties allow us to simplify lianas in the forests as we did in Section \ref{section:partitioned_methods}. Then we are left with forests with white nodes such as $\ptree 3 2 0 1 1$. This is where the last property comes in handy, as we can integrate by part this tree and obtain $\ptree 3 2 0 1 1 \sim -\ptree 3 1 0 1 3$.
We deduce the following theorem.

\begin{theorem}
\label{theorem:RK_partitioned_non_reversible}
Consider an ergodic Runge-Kutta method and a postprocessor as in Theorem \ref{theorem:RK_partitioned}, and suppose $f_1=-\nabla V$ and $f_2$ satisfies \eqref{equation:f2_perturbed}. Then under the notation and assumptions of Theorem \ref{Theorem:Postprocessed_Runge_Kutta_order_conditions}, if the coefficients of the method are chosen such that $\gamma \sim 0$ then the method has at least order $p+1$ for the invariant measure.
In particular, the conditions for order 1 and order 2 are in Table \ref{table:Order_conditions_for_partitioned_Runge-Kutta_method_with_postprocessor_and_perturbation}.

\begin{table}[tbh]
\setcellgapes{3pt}
\begin{center}
\begin{tabular}{|c|c|c|l|}
  \hline
  Order & Tree $\tau$ & $F(\tau)(\phi)$ & \multicolumn{1}{c|}{Order condition} \\
  \hline
  1 & $\etree 2 1 0 1$ & $\phi'f_1$ & $\sum b_i=1$ \\
  \hline
  2 & $\etree 3 1 0 1$ & $\phi'f_1'f_1$ & $\sum b_i c_i-2\sum b_i d_i-2\sum \overline{b_i}+2\overline{d_0}^2=-\frac{1}{2}$ \\
  \cline{2-4}
   & $\ptree 3 1 0 1 1$ & $\phi'f_1'f_2$ & $\sum b_i \widehat{c_i}+\sum \overline{b_i}-\sum \widehat{\overline{b_i}}=0$ \\
  \cline{2-4}
   & $\ptree 3 1 0 1 2$ & $\phi'f_2'f_1$ & $\sum \widehat{b_i} c_i-2\sum \widehat{b_i} d_i+\sum \widehat{b_i} -\sum \overline{b_i}-\sum \widehat{\overline{b_i}}+2\overline{d_0}^2=0$ \\
  \cline{2-4}
   & $\ptree 3 1 0 1 3$ & $\phi'f_2'f_2$ & $\sum \widehat{b_i} \widehat{c_i}-\frac{1}{2}\left(\sum \widehat{b_i}\right)^2=0$ \\
  \cline{2-4}
   & $\etree 2 1 1 1$ & $\phi'\Delta f_1$ & $\sum b_i d_i^2-2\sum b_i d_i-\sum \overline{b_i}+\overline{d_0}^2=-\frac{1}{2}$ \\
  \cline{2-4}
   & $\ptree 2 1 1 1 1$ & $\phi'\Delta f_2$ & $\sum \widehat{b_i} d_i^2-2\sum \widehat{b_i} d_i+\sum \widehat{b_i}-\sum \widehat{\overline{b_i}}+\overline{d_0}^2=0$ \\
  \hline
\end{tabular}
\caption{Order conditions for partitioned Runge-Kutta method with postprocessor for the perturbed equation (See Theorem \ref{theorem:RK_partitioned_non_reversible}). The sums are over all involved indices.}
\label{table:Order_conditions_for_partitioned_Runge-Kutta_method_with_postprocessor_and_perturbation}
\end{center}
\setcellgapes{1pt}
\end{table}
\end{theorem}

We note that if $f_2=J\nabla V$, we have $\frac{\sigma^2}{2} \ptree 2 1 1 1 1 \sim -\ptree 3 1 0 1 1 - \ptree 3 1 0 1 2$.
In this case, the order condition for $\ptree 2 1 1 1 1$ can be omitted and the two conditions of $\ptree 3 1 0 1 1$ and $\ptree 3 1 0 1 2$ are respectively replaced by
\begin{align*}
2\sum \widehat{b_i} d_i-\sum \widehat{b_i}+\sum b_i \widehat{c_i}-\sum \widehat{b_i} d_i^2+\sum \overline{b_i}-\overline{d_0}^2&=0,\\
\sum \widehat{b_i} c_i-\sum \widehat{b_i} d_i^2-\sum \overline{b_i}+\overline{d_0}^2&=0.
\end{align*}

\begin{remark}
In order for the method to satisfy $\AA_0=\LL$, the condition $\sum \widehat{b_i}=1$ should be added in Table \ref{table:Order_conditions_for_partitioned_Runge-Kutta_method_with_postprocessor_and_perturbation}, but it is not necessary to achieve order 1 for the invariant measure.
\end{remark}

\begin{ex*}
If $f_1$ satisfies Assumption \ref{Assumption:f_gradient} and $f_2$ satisfies \eqref{equation:f2_perturbed}, the following consistent postprocessed scheme has order 2 for the invariant measure (if it is ergodic):
\begin{equation}
\label{equation:exemple_non_reversible_perturbation}
\begin{array}{rl}
  X_{n+1}&=X_n+h f_1(X_n+\frac{\sigma}{2}\sqrt{h}\xi_n)+\frac{5}{4} h f_2(X_n+\frac{\sigma}{2}\sqrt{h}\xi_n)\\
  &-\frac{1}{4} h f_2(X_n-2 h f_2(X_n+\frac{\sigma}{2}\sqrt{h}\xi_n)-\frac{\sigma}{2}\sqrt{h}\xi_n)+\sigma\sqrt{h}\xi_n,\\
  \overline{X_n}&=X_n+\frac{\sigma}{2}\sqrt{h}\,\overline{\xi_n}.
\end{array}
\end{equation}
If $f_2=J\nabla V$, it needs two evaluations of $\nabla V$ per timestep similarly to a standard Runge-Kutta weak order 2 method.
For $f_2=0$, note that the scheme \eqref{equation:exemple_non_reversible_perturbation} coincides with the one proposed in \cite{Leimkuhler13rco}, formulated in a different manner (See \cite{Vilmart15pif}).
\end{ex*}


\bigskip

\noindent \textbf{Acknowledgements.}\
The authors would like to thank Hans Munthe-Kaas and Olivier Verdier for helpful discussions about an earlier version of this work.
This work was partially supported by the Swiss National Science Foundation, grants No. 200021\_162404, 200020\_178752 and 200020\_144313/1.

\bibliographystyle{abbrv}
\bibliography{Ma_Bibliographie}

\vskip-1ex

\appendix
\section*{Appendix}
\noindent Expression of $\widetilde{\gamma_2}$ in Theorem \ref{Theorem:Formal_IPP}.
\label{Appendix_gamma_2}
{\small
\begin{align*}
\widetilde{\gamma_2} &= \lpar a(\etree 4 1 0 1)-\frac{2}{\sigma^2}a(\etree 3 1 1 4)-\frac{2}{\sigma^2}a(\etree 3 2 1 2)+\frac{2}{\sigma^2}a(\etree 3 1 1 6)+\frac{2}{\sigma^2}a(\etree 3 2 1 3)-\frac{4}{\sigma^2}a(\etree 3 2 1 4)-\frac{4}{\sigma^4}a(\etree 2 1 2 6)\\[-0.5ex]
&+\frac{12}{\sigma^4}a(\etree 2 1 2 7)-\frac{24}{\sigma^6}a(\etree 1 1 3 1)\rpar \etree 4 1 0 1
+\lpar  a(\etree 3 1 1 1)-a(\etree 3 1 1 4)-a(\etree 3 2 1 2)+a(\etree 3 1 1 6)+a(\etree 3 2 1 3)-2 a(\etree 3 2 1 4)\\[-0.5ex]
&-\frac{2}{\sigma^2}a(\etree 2 1 2 6)+\frac{6}{\sigma^2}a(\etree 2 1 2 7)-\frac{12}{\sigma^4}a(\etree 1 1 3 1)\rpar \etree 3 1 1 1
+\lpar  a(\etree 4 2 0 1)-\frac{2}{\sigma^2} a(\etree 3 1 1 5)+\frac{4}{\sigma^4} a(\etree 2 1 2 2)+\frac{2}{\sigma^2}a(\etree 3 1 1 6)\\[-0.5ex]
&-\frac{4}{\sigma^4} a(\etree 2 1 2 6)+\frac{4}{\sigma^4} a(\etree 2 1 2 7)-\frac{8}{\sigma^6}a(\etree 1 1 3 1)\rpar \etree 4 2 0 1
+\lpar  a(\etree 3 1 1 2)-a(\etree 3 1 1 5)-a(\etree 3 2 1 2)-a(\etree 3 1 1 4)+2a(\etree 3 1 1 6)\\[-0.5ex]
&+\frac{2}{\sigma^2} a(\etree 2 1 2 2)+a(\etree 3 2 1 3)-2a(\etree 3 2 1 4)-\frac{4}{\sigma^2}a(\etree 2 1 2 6)+\frac{8}{\sigma^2} a(\etree 2 1 2 7)-\frac{16}{\sigma^4} a(\etree 1 1 3 1)\rpar \etree 3 1 1 2
+\lpar  a(\etree 3 1 1 3)\\[-0.5ex]
&-\frac{2}{\sigma^2}a(\etree 2 1 2 4)-a(\etree 3 1 1 5)+\frac{2}{\sigma^2}a(\etree 2 1 2 1)+a(\etree 3 1 1 6)+\frac{4}{\sigma^2} a(\etree 2 1 2 2)-\frac{4}{\sigma^2}a(\etree 2 1 2 6)+\frac{4}{\sigma^2} a(\etree 2 1 2 7)-\frac{8}{\sigma^4} a(\etree 1 1 3 1)\rpar \etree 3 1 1 3\\[-0.5ex]
&+\lpar  a(\etree 2 1 2 5)-a(\etree 2 1 2 4)+a(\etree 2 1 2 1)+ a(\etree 2 1 2 2)-a(\etree 2 1 2 6)+ a(\etree 2 1 2 7)-\frac{2}{\sigma^2} a(\etree 1 1 3 1)\rpar \etree 2 1 2 5
+\lpar  a(\etree 4 4 0 1)\\[-0.5ex]
&-\frac{2}{\sigma^2}a(\etree 3 2 1 4)+\frac{4}{\sigma^4} a(\etree 2 1 2 7)-\frac{8}{\sigma^6}a(\etree 1 1 3 1)\rpar \etree 4 4 0 1
+\lpar  a(\etree 3 2 1 1)-a(\etree 3 2 1 3)+2 a(\etree 3 2 1 4)-\frac{2}{\sigma^2}a(\etree 2 1 2 1)\\[-0.5ex]
&+\frac{2}{\sigma^2}a(\etree 2 1 2 6)-\frac{6}{\sigma^2}a(\etree 2 1 2 7)+\frac{12}{\sigma^4}a(\etree 1 1 3 1)\rpar \etree 3 2 1 1
+\lpar  a(\etree 4 3 0 1)-\frac{2}{\sigma^2}a(\etree 3 2 1 3)-\frac{2}{\sigma^2}a(\etree 3 1 1 6)+\frac{4}{\sigma^4}a(\etree 2 1 2 6)\\[-0.5ex]
&+\frac{4}{\sigma^2} a(\etree 3 2 1 4)-\frac{12}{\sigma^4}a(\etree 2 1 2 7)+\frac{24}{\sigma^6}a(\etree 1 1 3 1)\rpar \etree 4 3 0 1.
\end{align*}
}






\end{document}